\documentclass[12pt]{amsart}

\newtheorem{theorem}{Theorem}
\newtheorem{proposition}{Proposition}
\newtheorem{definition}{Definition}
\newtheorem{lemma}{Lemma}
\newtheorem{corollary}{Corollary}

\begin{document}

\title{The free unitary compact quantum group}

\author{Teodor Banica}
\address{Department of Mathematics, Paris 7 University, 75005 Paris, France}

\begin{abstract}
The free analogues of $U(n)$ in Woronowicz's compact quantum group theory are the quantum groups $\{A_u(F)|F\in GL(n,\mathbb C)\}$ introduced by Van Daele and Wang. We classify here their irreducible representations. Their fusion rules turn to be related to the combinatorics of Voiculescu's circular variable. If $F\bar{F}\in\mathbb R I_n$ we find an embedding $A_u(F)_{red}\subset C(\mathbb T)*_{red}A_o(F)$, where $A_o(F)$ is the deformation of $SU(2)$ that we previously studied. We use the representation theory and Powers' method for showing that the reduced algebras $A_u(F)_{red}$ are simple, with at most one trace.
\end{abstract}

\maketitle

\section*{Introduction}

One of the basic constructions in harmonic analysis is the Pontrjagin duality, which associates to any abelian group the abelian group of its characters. This duality has been extended to the non-abelian groups, but the dual object (the convolution algebra of the group) is no longer of the same nature. In order to have a framework covering at the same time the groups and their duals, we are led to the construction of new objects in the category of Hopf algebras, called ``quantum groups''.

A certain number of families of examples have been studied at the operator algebra level, and in particular, Woronowicz \cite{wo2} constructed in 1987 the class of ``compact matrix quantum groups''. Such a quantum group is a pair $(G,u)$ formed by a unital $C^*$-algebra $G$ and a matrix $u\in M_n(G)$ such that:
\begin{enumerate}
\item The coefficients of $u$ generate a dense $*$-subalgebra $G_s\subset G$.

\item There exists a $C^*$-morphism $\delta:G\to G\otimes_{min}G$ mapping  $u_{ij}\to\sum_ku_{ik}\otimes u_{kj}$.

\item There exists a linear antimultiplicative map $\kappa:G_s\to G_s$ such that $\kappa(\kappa(a^*)^*)=a$ for any $a\in G_s$, and such that $(Id\otimes\kappa)u=u^{-1}$.
\end{enumerate}

This definition covers as well the ``compact quantum'' case, obtained by taking projective limits, and the ``discrete quantum'' case too, obtained by duality. The general ``locally compact quantum'' case was treated by Baaj and Skandalis in \cite{bsk}.

Given $n\in\mathbb N$, the universal $C^*$-algebra $A_u(I_n)$ generated by the coefficients of a unitary $n\times n$ matrix, whose transpose is also unitary, is a compact matrix quantum group \cite{vwa}, \cite{wa1}, \cite{wa2}. $A_u(I_n)$ is an analogue of $U(n)$ in Woronowicz's theory. This algebra, and its ``deformed'' versions $\{A_u(F)|F\in GL(n,\mathbb C)\}$, are the objects to be investigated here.

\medskip

{\bf Acknowledgements.} I would like to express my deepest gratitude to my PhD advisor, G. Skandalis. I would like to thank as well E. Blanchard for many useful discussions on the Hopf $C^*$-algebras, and S. Wang for several useful comments on this paper.

\section{Definitions, and statements of the results}

In this section we construct the compact matrix quantum groups $A_u(F)$, in a slightly different manner than in the article of Van Daele and Wang \cite{vwa}, and we state our main results. The end of this section contains the plan of the article, reminders and notations.

\medskip

(1) There are several possible definitions for the morphisms between the compact matrix quantum groups, corresponding to different notions of isomorphism. Without getting into details (in this paper we write $(G,u)=(H,v)$ when $G=H$ as $C^*$-algebras, and when $u=v$) let us recall from \cite{wo2} the definition of the similarity:

\medskip

{\em Two compact matrix quantum groups $(G,u)$ and $(H,v)$, with $u\in M_n(G),v\in M_m(H)$, are called similar (and we write $G\sim_{sim}H$) if $n=m$, and there exists a matrix $Q\in GL(n,\mathbb C)$ and a $C^*$-isomorphism $f:G\to H$ such that $(Id\otimes f)u=QvQ^{-1}$.}

\medskip

(2) Let $(G,u)$ be a compact matrix quantum group. We call representation of $(G,u)$ any invertible matrix $r\in M_k(G)$ satisfying:
$$(Id\otimes\delta)r=r_{12}r_{13}:=\sum e_{ij}\otimes r_{ik}\otimes r_{kj}$$

Woronowicz's Peter-Weyl type theory in \cite{wo2} shows that any representation is equivalent to a unitary representation. In particular, $v=Q^{-1}uQ$ is unitary, for a certain matrix $Q\in GL(n,\mathbb C)$. Up to replacing $(G,u)$ by a similar compact matrix quantum group, we can assume that $u$ is unitary.

\medskip

(3) Let $(G,u)$ be a compact matrix quantum group, with $u\in M_n(G)$ unitary. The representation $\bar{u}=(u_{ij}^*)$ is then equivalent to a unitary representation, so there exists $F\in GL(n,\mathbb C)$ such that $F\bar{u}F^{-1}$ is unitary. Thus $G$ is a quotient of $A_u(F)$, where:

\begin{definition}
To any $n\in\mathbb N$ and any $F\in GL(n,\mathbb C)$ we associate the universal $C^*$-algebra $A_u(F)$ generated by variables $\{u_{ij}\}_{i,j=1,\ldots,n}$, with the relations making unitaries the matrices $u=(u_{ij})$ and $F\bar{u}F^{-1}$.
\end{definition}

Observe that $A_u(F)$ is well-defined. Indeed, if $J$ is the two-sided ideal generated inside the free algebra on $2n^2$ variables $L=\mathbb C<u_{ij},u_{ij}^*>$ by the relations making unitaries the matrices $u=(u_{ij})$ and $F\bar{u}F^{-1}=F(u_{ij}^*)F^{-1}$, then the images of the generators $u_{ij},u_{ij}^*$ in the quotient $L/J$ have norm $\leq1$, for any $C^*$-norm on $L/J$. Thus, $L/J$ has an enveloping $C^*$-algebra, that we can denote $A_u(F)$.

$(A_u(F),u)$ is a compact matrix quantum group. Indeed, we have $v$ unitary $\implies$ $v_{12}v_{13}$ unitary, and by applying this to $v=u$ and to $v=F\bar{u}F^{-1}$ (with the remark that we have $F\overline{u_{12}u_{13}}F^{-1}=(F\bar{u}F^{-1})_{12}(F\bar{u}F^{-1})_{13}$) we can define $\delta$ by using the universal property. As for the existence of the antipode $\kappa$, this follows from Woronowicz's result in \cite{wo4}, stating that this condition is equivalent to the fact that both $u,\bar{u}$ are invertible.

\medskip

Observe that for any compact group $G\subset U(n)$, the algebra $C(G)$ is a quotient of $C(U(n))$. By the above, the analogue of $U(n)$ among the compact quantum groups is the {\em family} $\{A_u(F)|F\in GL(n,\mathbb C)\}$.

\medskip

{\bf Remark.} The relations which define $A_u(F)$ are as follows:
$$uu^*=u^*u=(F^*F)\bar{u}(F^*F)^{-1}u^t=u^t(F^*F)\bar{u}(F^*F)^{-1}=I$$ 

We deduce from this that we have the following equalities:
$$A_u(F)=A_u(\sqrt{F^*F})=A_u(\lambda F),\quad\forall F\in GL(n,\mathbb C),\forall\lambda\in\mathbb C^*$$

There are as well a number of similarities between $A_u(F)$'s, for instance for $V,W\in U(n)$ and $F\in GL(n,\mathbb C)$ we have $A_u(F)\sim_{sim}A_u(VFW)$ (see Proposition 6 below). Thus, we could have used other parameters for $A_u(F)$, such as $F^*F$, or $\sqrt{F^*F}$, or the eigenvalue list of $\sqrt{F^*F}$ and so on, see \cite{vwa}, \cite{wa2}. Of course, the choice of the parameter is not a serious issue, because we always obtain the same objects, at least up to similarity.

\medskip

The quotient of $A_u(F)$ by the relations $u=F\bar{u}F^{-1}$ can be thought of as being an ``orthogonal version'' of $A_u(F)$. Observe that the condition $u=F\bar{u}F^{-1}$ implies $\bar{u}=\bar{F}u\bar{F}^{-1}$, and so $u=(F\bar{F})u(F\bar{F})^{-1}$. Thus, if $F\bar{F}$ is not a scalar multiple of the identity of $M_n(\mathbb C)$, then $u$ is not irreducible in this quotient. Observe also that $F\bar{F}=cI_n$ with $c\in\mathbb C$ implies $\bar{F}F=\bar{c}I_n$, and so $c=\bar{c}$. These observations lead to:

\begin{definition}
To any $n\in\mathbb N$ and any $F\in GL(n,\mathbb C)$ satisfying $F\bar{F}=cI_n$ with $c\in\mathbb R$ we associate the quotient $A_o(F)$ of $A_u(F)$ by the relations $u=F\bar{u}F^{-1}$.
\end{definition}

The irreducible representations of $A_o(F)$ are known to be indexed by $\mathbb N$, and their fusion rules coincides with the fusion rules for $SU(2)$ (cf. \cite{ba2}, see also Theorem 4 below).

\medskip

{\bf Notations.} $\mathbb N*\mathbb N$ is the coproduct in the category of monoids of two copies of $\mathbb N$, generated respectively by $\alpha,\beta$. We denote by $e$ the unit of $\mathbb N*\mathbb N$, and by $-$ the involution of $\mathbb N*\mathbb N$ given by $\bar{e}=e,\bar{\alpha}=\beta,\bar{\beta}=\alpha$ and antimultiplicativity.

\medskip

The main result in this paper is as follows:

\begin{theorem}
Let $n\in\mathbb N$, and let $F\in GL(n,\mathbb C)$.
\begin{enumerate}
\item The irreducible representations of $(A_u(F),u)$ can be labelled by $\mathbb N*\mathbb N$, with $r_e=1$, $r_\alpha=u$, $r_\beta=\bar{u}$. For any $x,y\in\mathbb N*\mathbb N$ we have $r_{\bar{x}}=\bar{r}_x$ and:
$$r_x\otimes r_y=\sum_{\{a,b,g\in\mathbb N*\mathbb N|x=ag,y=\bar{g}b\}}r_{ab}$$

\item The subalgebra of $A_u(F)$ generated by the characters of the representations is the free $*$-algebra on the character $\chi(u)$ of the fundamental representation.

\item $\chi(u)/2$ is a circular variable over $A_u(F)$, with respect to the Haar measure.

\item If $F\bar{F}\in\mathbb RI_n$ we have $A_u(F)_{red}\subset C(\mathbb T)*_{red}A_o(F)$ via $u_{ij}\to zv_{ij}$, where $v$ is the fundamental representation of $A_o(F)$, and $z$ is the standard generator of $C(\mathbb T)$.
\end{enumerate}
\end{theorem} 

Observe that (1) above shows that the family $\mathcal F=\{A_u(F)|n\in\mathbb N,F\in GL(n,\mathbb C)\}$ has the following remarkable property:

\medskip

{\em For any $G,H\in\mathcal F$ there exists a bijection $\psi$ between the equivalence classes of representations of $G$, and those of $H$, which preserves the sums and the tensor products, and which sends the irreducible representations of $G$ to those of $H$, and the fundamental representation of $G$ to that of $H$.}

\medskip

An important result of this type, for the family with parameter $q>0$ of compact quantum groups associated to a classical Lie algebra, was proved by Rosso \cite{ro1}, \cite{ro2}. Another result of this type, this time of ``rigidity'', is the one from \cite{ba2}, stating that the following family has the above property, and is in addition {\em maximal} with this property:
$$\left\{A_o(F)\big|n\in\mathbb N,F\in GL(n,\mathbb C),F\bar{F}\in\mathbb RI_n\right\}$$

The following result, that we will prove here, is of the same type:

\begin{theorem}
Assuming that the irreducible representations of a compact matrix quantum group $(G,u)$ can be labelled by $\mathbb N*\mathbb N$, with $r_e=1$, $r_\alpha=u$, $r_\beta=\bar{u}$ and $r_x\otimes r_y=\sum_{x=ag,y=\bar{g}b}r_{ab}$, there exist $n\in\mathbb N$ and $F\in GL(n,\mathbb C)$ such that $G_{full}\sim_{sim}A_u(F)$.
\end{theorem}

The representation theory of $A_u(F)$ is the subject of the first part (sections 2,3,4) of this paper. In the second part (sections 6,7,8) we use representation theory for solving certain topological questions regarding the algebras $A_u(F)$ and $A_u(F)_{red}$.

Let us recall that for a compact matrix quantum group $(G,u)$ the Haar functional is not necessarily a trace, but satisfies a formula of the following type:
$$h(xy)=h(y(f_1*x*f_1)),\quad\forall x,y\in G_s$$

To be more precise, here $*$ is the convolution over $G_s$, and $\{f_z\}_{z\in\mathbb C}$ is a certain family of characters of $G_s$ (see Theorem 5.6 in \cite{wo2}). Our third and last main result will be:

\begin{theorem}
Let $n\in\mathbb N$ and $F\in GL(n,\mathbb C)$. Then $A_u(F)_{red}$ is simple.

Let $s,t\in\mathbb R$ and and let $\psi$ be a state of $A_u(F)_{red}$ such that $\psi(xy)=\psi(y(f_s*x*f_t))$, for any $x,y\in A_u(F)_{red}$. Then $\psi$ must be the Haar measure of $A_u(F)_{red}$.

In particular, if $F$ is the scalar multiple of a unitary matrix, then $A_u(F)_{red}$ is simple with unique trace. Otherwise, $A_u(F)_{red}$ is simple without trace.
\end{theorem}

Among the other results regarding $A_o(F)$ and $A_u(F)$, let us mention:
\begin{enumerate}
\item A commutation result inside $A_u(I_2)$.

\item The equality of factors $A_u(I_2)_{red}''=W^*(F_2)$.

\item The non-amenability of $A_o(F)$ and $A_u(F)$.

\item The non-nuclearity of $A_o(I_n)_{red}$ and $A_u(I_n)_{red}$.

\item Some remarks on the Haar measures of $A_o(F)$ and $A_u(F)$.
\end{enumerate}

Some of these results are in fact particular cases of certain more general statements, regarding the compact quantum groups. Let us mention here the simplicity result (Proposition 9 below), which in the case of $G=C^*_{red}(F_n)$ contains a simplification with respect to the known proofs of the simplicity of $C^*_{red}(F_n)$, from \cite{har}, \cite{hsk}, \cite{pow}.

\medskip

The paper is organized as follows:

-- $2^{nd}$ section: we recall here our results in \cite{ba2} regarding $A_o(F)$, and we give a description (in terms of certain noncrossing partitions) of the space of fixed vectors for $u^{\otimes k}$.

-- $3^{rd}$ section: we construct here the abstract algebra $A$ generated by variables $\{r_x|x\in\mathbb N*\mathbb N\}$ which multiply according to the formulae $r_xr_y=\sum_{x=ag,y=\bar{g}b}r_{ab}$, and we show that we have $A\simeq\mathbb C<X,X^*>$. By using this algebra, and the same method as in the ``orthogonal'' case, from \cite{ba2}, we conclude that the proof of Theorem 1 is equivalent to the computation of the dimensions of the commutants of the following representations:
$$u^{\otimes m_1}\otimes\bar{u}^{\otimes n_1}\otimes u^{\otimes m_2}\otimes\bar{u}^{\otimes n_2}\otimes\ldots\qquad (\star)$$

These dimensions are $*$-moments of the character $\chi(u)$ of the fundamental representation, with respect to the Haar integration, and we see that $\chi(u)/2$ must be {\em circular}.

-- $4^{th}$ section: in the case $F\bar{F}\in\mathbb RI_n$, we combine the results on $A_o(F)$ with a free probability result, for proving Theorem 1. Then, we use the reconstruction results in \cite{wo3} in order to find a system of generators for the spaces of fixed vectors of representations of type ($\star$). The dimensions of these spaces are exactly the $*$-moments of $\chi(u)$, and by using this remark, we can pass from the case $F\bar{F}\in\mathbb RI_n$ to the general case, $F\in GL(n,\mathbb C)$.

-- $5^{th}$ section: we describe here $A_o(F)$ and $A_u(F)$ in the case $F\in GL(2,\mathbb C)$. Theorem 1 (4) allows us to identify $A_u(I_2)_{red}$ as a $C^*$-subalgebra of $C(\mathbb T)*_{red}C(SU(2))$, and we deduce from this two Hopf $C^*$-algebra embeddings $C(SO(3))\subset A_u(I_2)$, as well as the equality of factors $A_u(I_2)_{red}''=W^*(F_2)$.

-- $6^{th}$ section: here we use the characters $\{f_z\}$ from \cite{wo2} in order to ``deform'' the adjoint representation of a compact matrix quantum group.

-- $7^{th}$ section: we generalize here to the case of the compact quantum groups the ``Powers' property'' of de la Harpe \cite{har}, and the proof of simplicity from \cite{hsk}. The idea is that of replacing the inner automorphisms $x\to u_gxu_g^*$ from the discrete case by completely positive maps of type $x\to ar(r)x$, with $r\in\widehat{A}$.

-- $8^{th}$ section: $A_u(F)_{red}$ does not have Powers' property, but by using the computations from sections 6 and 7, we manage however to prove Theorem 3.

\medskip

{\bf Notations and Reminders:}

\medskip

{\em (A) Matrices.} We denote by $\{e_1,\ldots,e_n\}$ the standard basis of $\mathbb C^n$, and by $e_{ij}$ the standard matrix units, given by $e_{ij}:e_j\to e_i$. If $A$ is a $*$-algebra and $u\in M_n(A)$, $u=\sum e_{ij}\otimes u_{ij}$, we set $\bar{u}=\sum e_{ij}\otimes u_{ij}^*$, $u^t=\sum e_{ij}\otimes u_{ji}$, $u^*=\sum e_{ij}\otimes u_{ji}^*$.

\medskip

{\em (B) Representations.} Given a compact quantum group $G$, we denote by $Rep(G)$ the set of equivalence classes of unitary representations of $G$, and by $\widehat{G}\subset Rep(G)$ the set of equivalence classes of unitary irreducible representations.

If $u=\sum e_{ij}\otimes u_{ij}\in M_n(G)$ and $v=\sum e_{ij}\otimes v_{ij}\in M_m(G)$ are two representations, we denote by $u\otimes v$ the matrix $u_{13}v_{23}=\sum e_{ij}\otimes e_{kl}\otimes u_{ij}v_{kl}$, and by $u+v$ the matrix $diag(u,v)$, which are both representations. The map $(u,v)\to u\otimes v$ makes $Rep(G)$ a monoid. The same holds for $(u,v)\to u+v$.

The character of $u$ is $\chi(u)=\sum u_{ii}\in G$. We have $\chi(u+v)=\chi(u)+\chi(v)$ and $\chi(u\otimes v)=\chi(u)\chi(v)$. See \cite{wo2}, \cite{wo3}.

\medskip

{\em (C) Woronowicz's Peter-Weyl theory.} We denote by $G_{central}$ the linear space (and so, $*$-algebra) generated inside the ``coefficient'' $*$-algebra $G_s$ by the characters of the representations. We will often use, without reference, the following fundamental result (Theorem 5.8 in \cite{wo2}):

{\em The Haar measure is a trace on $G_{central}$. The set $\{\chi(u)|u\in\widehat{G}\}$ is a basis of $G_{central}$, which is orthonormal with respect to the scalar product associated to the Haar measure.}

\medskip

{\em (D) Full and reduced version.} The reduced version of a compact matrix quantum group $(G,u)$ is $G_{red}=G/\{x|h(xx^*)=0\}$, with $h$ being the Haar measure of $G$. The full version is $G_{full}=C^*(G_s)$, the enveloping $C^*$-algebra of $G_s$. Both $G_{full},G_{red}$ are compact matrix quantum groups (cf. \cite{bsk}, \cite{wo2}). $G$ is said to be amenable when the projection $G_{full}\to G_{red}$ is an isomorphism, and is called full (resp. reduced) when the projection $G_{full}\to G$ (resp. $G\to G_{red}$) is an isomorphism. We have $G_s=H_s\iff G_{red}=H_{red}\iff G_{full}=H_{full}$. Observe also that $A_u(F)$ and $A_o(F)$ are both full, by definition.

\medskip

{\em (E) Freeness.} If $(A,\phi)$ is a unital $*$-algebra endowed with a unital linear form, a family of subalgebras $1\in A_i\subset A$ ($i\in I$) is called free when $a_j\in A_j\cap\ker(\phi)$ with $i_j\neq i_{j+1},\, 1\leq j\leq n-1$ implies  $a_1a_2\ldots a_n\in\ker(\phi)$ (see \cite{vdn}). Two elements $a,b\in A$ are called $*$-free when the unital $*$-algebras that they generate inside $A$ are free. Fundamental example: let $(A,\phi)$ and $(B,\psi)$ be two unital $C^*$-algebras with states, and consider their free product $A*B$ (coproduct in the category of unital $C^*$-algebras). If we denote by $\phi*\psi$ the free product of $\phi,\psi$ and by $\pi_{\phi*\psi}$ the associated GNS representation, then $\pi_{\phi*\psi}(A)$ and $\pi_{\phi*\psi}(B)$ are free inside $\pi_{\phi*\psi}(A*B)$ (see \cite{avi}, \cite{voi}, \cite{vdn}).

\medskip

{\em (F) Free products.} If $(G,u)$ and $(H,v)$ are two compact matrix quantum groups, then $(G*H,diag(u,v))$ is a compact matrix quantum group too, and its Haar measure is the free product $h*k$ of the Haar measures of $G,H$ (see \cite{wa2}). The reduced free product $\pi_{h*k}(G*H)$, which is a reduced compact matrix quantum group, will be denoted $G*_{red}H$. Since $h,k$ are faithful on $G_{red},H_{red}$, we have embeddings $G_{red},H_{red}\subset G*_{red}H$.

\medskip

{\em (G) $*$-distributions.} The $*$-distribution of an element $a\in (M,\phi)$ of a $*$-algebra with a linear form is the functional on $\mathbb C<X,X^*>$ given by $P\to\phi(P(a,a^*))$:
$$\mu_a=\left[\mathbb C<X,X^*>\,\, \displaystyle{\mathop{\longrightarrow}^{X\mapsto a}}\,\, M\,\, \displaystyle{\mathop{\longrightarrow}^\phi}\,\, \mathbb C\right]$$

The $*$-moments of $a$ are the values of $\mu_a$ on the noncommutative monomials in $X,X^*$, i.e. on the monoid generated inside $\mathbb C<X,X^*>$ by the elements $X,X^*$.

If $(M,\phi)$ is a $C^*$-algebra with a faithful state, and $a=a^*$, then $\mu_a$ can be viewed (by restricting to $\mathbb C[X]$, then completing) as a probability measure on the spectrum of $a$.

\medskip

{\em (H) Circular variables.} The centered semicircle law is the measure $\gamma_{0,1}=\frac{2}{\pi}\sqrt{1-t^2}dt$ on $[-1,1]$. Any self-adjoint element having this distribution is called semicircular. A quart-circular element is a positive element having as distribution the measure $\frac{4}{\pi}\sqrt{1-t^2}dt$ on $[0,1]$. A Haar unitary is a unitary $u$ such that $\mu_u(X^k)=0$ for any $k\neq 0$. A variable $g$ is called circular if $2^{-1/2}(g+g^*)$ and $-i2^{-1/2}(g-g^*)$ are semicircular and free (see \cite{vdn}).

\section{More about $A_o(F)$}

It is easy to see, starting from the definition of $A_o(F)$, that we have $A_0\begin{pmatrix}0&1\\-1&0\end{pmatrix}=C(SU(2))$. Moreover, modulo similarity, we have an equality as follows:
$$\left\{A_o(F)\big|F\in GL(2,\mathbb C),F\bar{F}\in\mathbb RI_2\right\}=\left\{S_\mu U(2)\big|\mu\in[-1,1]-\{0\}\right\}$$

Here $S_\mu U(2)$ are the deformations of $S_1U(2)=C(SU(2))$ constructed by Woronowicz in \cite{wo1}, \cite{wo2} (see section 5 below). For $F\in GL(n,\mathbb C)$ with $n$ arbitary, we have:

\begin{theorem}[\cite{ba2}]
Let $n\in\mathbb N$ and $F\in GL(n,\mathbb C)$ satisfying $F\bar{F}\in\mathbb RI_n$. Then the irreducible representations of $A_o(F)$ are self-adjoint, and can be indexed by $\mathbb N$, with $r_0=1$, $r_1=u$ and
$$r_kr_s=r_{|k-s|}+r_{|k-s|+2}+\ldots+r_{k+s-2}+r_{k+s}$$
(i.e. the same formulae as for the representations of $SU(2)$).
\end{theorem}

Let us briefly recall the proof. The condition $F\bar{F}\in\mathbb RI_n$ shows that the projection onto $\mathbb C\sum F_{ji}e_i\otimes e_j$, which intertwines $u^{\otimes 2}$, defines for any $k$ a representation of the Jones algebra $A_{\beta,k}$ \cite{jon} into $Mor(u^{\otimes k},u^{\otimes k})$. By using  \cite{wo3} one can prove that this representation is surjective, and the inequality $\dim(Mor(u^{\otimes k},u^{\otimes k}))\leq\dim(A_{\beta,k})\leq C_k$ that we obtain in this way allows us to construct (by recurrence on $k$) irreducible representations $r_k$ of $A_o(F)$ which satisfy the same multiplication formulae as those of $SU(2)$.

As a corollary of the proof (see Remark (2) in \cite{ba2}), we have:
$$\dim(Mor(u^{\otimes k},u^{\otimes k}))=\dim(A_{\beta ,k})=C_k=\frac{(2k)!}{k!(k+1)!}$$

Let $h$ be the Haar measure of $A_o(F)$. By Reminder (C), we have:
$$\dim(Mor(u^{\otimes k},u^{\otimes k}))=h(\chi(u)^{\otimes 2k})=\dim(Mor(1,u^{\otimes 2k}))$$

The Catalan numbers $C_k$ have another remarkable property: they are the moments of the semicircle law $\gamma_{0,1}$ of Wigner and Voiculescu. Indeed, one can compute the moments of $\gamma_{0,1}$ by using 3.3 and 3.4 in \cite{vdn}, and the formula of the residues, as follows:
$$\gamma_{0,1}(X^{2k})=(2k+1)^{-1}(2\pi i)^{-1}\int_T (z^{-1}+z/4)^{2k+1}=4^{-k}\frac{(2k)!}{k!(k+1)!}$$ 

By combining these equalities, we obtain:

\begin{proposition}
$\chi(u)/2\in(A_o(F),h)$ is a semicircular variable.
\end{proposition}

{\bf Remark.} In the case $F=\begin{pmatrix}0&1\\-1&0\end{pmatrix}$, the character of the fundamental representation $u=\begin{pmatrix}a&b\\-\bar{b}&\bar{a}\end{pmatrix}$ of $A_o(F)=C(SU(2))$ is $\chi(u)=2Re(a)$. The fact that $Re(a)$ is semicircular with respect to the Haar measure of $SU(2)$ can be viewed geometrically, by identifying $SU(2)$ with the sphere $S^3$, and its Haar measure with the uniform measure on this sphere.

\medskip

The above result has the following consequence:

\begin{corollary}[G. Skandalis]
$A_o(F)$ is amenable for $F\in GL(2,\mathbb C)$, and non-amenable for  $F\in GL(n,\mathbb C)$ with $n>2$.
\end{corollary}

\begin{proof}
The support of $\gamma_{0,1}$ being $[-1,1]$, and $h$ being faithful on $A_o(F)_{red}$, we obtain $Spec(\chi(u)/2)=[-1,1]$ inside $A_o(F)_{red}$ (see the Reminders). Thus, for $n\geq3$, the element $n-\chi(u)$ is invertible inside $A_o(F)_{red}$. But the counit of $A_o(F)$ is a unital $*$-morphism mapping $n-\chi(u)\to0$, and so $A_o(F)\neq A_o(F)_{red}$.

Finally, in the case $F\in GL(2,\mathbb C)$, we have a similarity between $A_o(F)$ and a certain $S_\mu U(2)$ (see section 5 below), which is amenable, cf. \cite{bla}, \cite{nag}.
\end{proof}

{\bf Remark.} A part of the classical results regarding the amenability was extended to the locally compact quantum group case in \cite{bsk}, \cite{bla} (see also Proposition 10 below). The above proof of the non-amenability of $A_o(F)$ can be extended to arbitrary compact matrix quantum groups, the result being that $(G,u)$ with $u\in M_n(G)$ is amenable if and only if the support of the law of $Re(\chi(u))$ with respect to the Haar mesure contains $n$.

\medskip

We recall that for any representation $r\in B(H_r)\otimes G$ of a compact quantum group $G$, the fixed vectors of $r$ are the elements $x\in H_r$ satisfying $r(x\otimes1)=x\otimes1$. These vectors form a linear subspace of $H_r$, which can be identified with $Mor(1,r)$.

We will describe now the fixed vectors for the representation $u^{\otimes k}$ of $A_o(F)$. The results which follow will be used in section 4 below, for proving Theorem 1 for the matrices $F$ which do not (!) verify the condition $F\bar{F}\in\mathbb RI_n$. So, the reader interested only in the algebras $A_u(F)$ with $F\bar{F}\in\mathbb RI_n$ (e.g. $A_u(I_n)$) can skip the rest of this section.

In order to describe the above-mentioned fixed vectors, we will need:

\begin{lemma}[\cite{ba2}] 
Let $n\in\mathbb N$ and $F\in GL(n,\mathbb C)$ satisfying $F\bar{F}\in\mathbb RI_n$. Let also $H=\mathbb C^n$, with orthonormal basis $\{e_i\}$. 
\begin{enumerate}
\item The operator $E\in B(\mathbb C,H^{\otimes 2})$ mapping $x\to x\sum e_i\otimes Fe_i$ belongs to $Mor(1,u^{\otimes 2})$.

\item We have $(E^*\otimes Id_H)(Id_H\otimes E)\in\mathbb CId_H$.

\item Let $Mor(r,s)\subset B(H^{\otimes r},H^{\otimes s})$ be the sets of linear combinations of (composable) products of maps of the form $Id_{H^{\otimes k}}$ or $Id_{H^{\otimes k}}\otimes E\otimes
Id_{H^{\otimes p}}$ or $Id_{H^{\otimes k}}\otimes E^*\otimes Id_{H^{\otimes p}}$. Then the concrete monoidal $W^*$-category of representations of $A_o(F)$ is the completion in the sense of \cite{wo3} of the following concrete monoidal $W^*$-category:
$$W(F):=\{\mathbb N,+,\{ H^{\otimes r}\}_{r\in\mathbb N},\{ Mor(r,s)\}_{r,s\in\mathbb N}\}$$
\end{enumerate}
\end{lemma}

By keeping the same notations, we have as well:

\begin{lemma}
Let $I(p)=Id_{H^{\otimes p}}$ and $V(p,q)=I(p)\otimes E\otimes I(q)$. 

\begin{enumerate}
\item Any morphism of $W(F)$ is a linear combination of maps of the form $I(.)$ or of the form $V(.,.)\circ\ldots\circ V(.,.)\circ V(.,.)^*\circ\ldots\circ V(.,.)^*$.

\item For any $k\geq 0$, the maps of the form $M\otimes I(1)\otimes N$ with $M\in Mor(0,2x)$, $N\in Mor(0,2y)$ and $x+y=k$ generate $Mor(1,2k+1)$.

\item For any $k\geq 0$, the maps of the form $(I(1)\otimes M\otimes I(1)\otimes N)\circ E$ with $M\in Mor(0,2x)$, $N\in Mor(0,2y)$ and $x+y=k$ generate $Mor(0,2k+2)$.
\end{enumerate}
\end{lemma}

\begin{proof}
Here (1), i.e. the fact that we can ``move the $*$ signs to the right'' follows from Lemma 1 (2). We prove now (2), by recurrence on $k$. At $k=0$ this is clear from (1), which gives $Mor(1,1)=\{\mathbb C I(1)\}$. So, let $k\geq1$ and $A\in Mor(1,2k+1)$. By (1), $A$ is a linear combination of maps of the form $V(k_1,s_1)\circ\ldots\circ V(k_m,s_m)$, with $k_1+s_1=2k-1$ and $T:=V(k_2,s_2)\circ\ldots\circ V(k_m,s_m)$, from $Mor(1,2k-1)$. By recurrence $T$ is a linear combination of maps of the form $B=(M\otimes I(1)\otimes N)$ for certain $M\in Mor(0,2x)$ and $N\in Mor(0,2y)$, with $x+y=k-1$, and we can conclude, as follows: 

-- either $k_1\geq 2x+1$, so $V(k_1,s_1)\circ B=M\otimes I(1)\otimes ((I(k_1-2x-1)\otimes E\otimes I(s_1))\circ N)$.

-- or $k_1\leq 2x$, and so $V(k_1,s_1)\circ B=((I(k_1)\otimes E\otimes I(2x-k_1))\circ M)\otimes I(1)\otimes N$.

We prove now (3). Let $A\in Mor(0,2k+2)$. By (1), $A$ is a linear combination of maps of the form $B=V(k_1,s_1)\circ\ldots\circ V(k_m,s_m)$. Since $V(k_m,s_m)=V(0,0)$, we can consider the smallest $p$ such that $k_p=0$. Then $B=(I(1)\otimes G)\circ(E\otimes I(s_p))\circ K$, with $G=V(k_1-1,s_1)\circ\ldots\circ V(k_{p-1}-1,s_{p-1})$
and $K=V(k_{p+1},s_{p+1})\circ\ldots\circ V(k_m,s_m)$. Thus:
$$B=(I(1)\otimes G)\circ (I(2)\otimes K)\circ E=(I(1)\otimes (G\circ (I(1)\otimes K)))\circ E$$

But $G\circ (I(1)\otimes K)\in Mor(1,2k+1)$ is, by (2), of the form $M\otimes I(1)\otimes N$, for certain $M\in Mor(0,2x)$ and $N\in Mor(0,2y)$, and so (3) holds indeed.
\end{proof}

We can now construct the basis that we are looking for, as follows:

\begin{proposition}
We define $W_{2k}(F)\subset Mor(0,2k)$ by $W_0(F)=1$ and (by recurrence) by:
$$W_{2k+2}(F)=\bigcup_{k=x+y}\{ (I(1)\otimes M\otimes I(1)\otimes N)\circ E\mid M\in W_{2x}(F), N\in W_{2y}(F)\}$$

Then $W_{2k}(F)$ is a basis of $Mor(0,2k)$, $\forall\, k\geq 0$. 
\end{proposition}

\begin{proof}
The numbers $D_k=Card(W_{2k}(F))$ verify $D_0=D_1=1$ and:
$$D_{k+1}=\sum_{k=x+y}D_xD_y$$

It follows that these are the Catalan numbers (well-known: consider the square of the series $\sum D_kz^k$\ldots). Thus $Card(W_{2k}(F))=dim(Mor(0,2k))$. On the other hand, Lemma 2 (3) shows that $W_{2k}(F)$ generates $Mor(0,2k)$, and this gives the result.
\end{proof}

{\bf Remark.} Let $P=P_1\coprod\ldots\coprod P_k$ be a noncrossing pairing of $\{1,\ldots,2k\}$, i.e. a pairing such that, with $P_m=\{ i_m,j_m\}$ with 
$i_m<j_m$ for any $1\leq m\leq k$, we have:
$$\forall\, m\neq n,\, i_m<i_n<j_m\,\implies j_n<j_m$$

We can associate to $P$ the following vector of $(\mathbb C^n)^{\otimes 2k}$ :
$$v(P)=\sum_{1\leq s_1\ldots s_{2k}\leq k}F_{s_{j_1}s_{i_1}}\ldots F_{s_{j_k}s_{i_k}}e_{s_1}\otimes\ldots \otimes e_{s_{2k}}$$

One can show by recurrence on $k$, by using Proposition 2, that the set of these vectors $v(P)$ coincides with the set $\{ X(1)|X\in W_{2k}(F)\}$, and so is a basis of the space of fixed vectors of the representation $u^{\otimes 2k}$ of $A_o(F)$. This allows in principle to compute the Haar measure of $A_o(F)$, because for any representation $r$ of a compact matrix quantum group we have $(Id\otimes h)(r)=$ projection on the space of fixed vectors of $r$ (see \cite{wo2}).

\section{Reconstruction of $A_u(F)_{central}$}

We construct and study here the algebra $A$ generated by variables $\{r_x|
x\in\mathbb N*\mathbb N\}$ which multiply according to the formula $r_xr_y=\sum_{x=ag,y=\bar{g}b}r_{ab}$.

\medskip

{\bf Notations. } $\mathbb N*\mathbb N$ is the free product (coproduct in the category of monoids) of two copies of $\mathbb N$, denoted multiplicatively $\{e,\alpha,\alpha^2,\ldots\}$ and $\{e,\beta,\beta^2,\ldots\}$. We consider the set $A$ of functions $\mathbb N*\mathbb N\to\mathbb C$ having finite support. We will identify $\mathbb N*\mathbb N\subset A$, as Dirac masses. With the addition and the multiplication of the functions, $A$ is the algebra of noncommutative polynomials in 2 variables, i.e. we have an isomorphism, as follows:
$$(\mathbb C<X,X^*>,+,\cdot )\simeq(A,+,\cdot )\quad,\quad X\to\alpha\, ,\,
X^*\to\beta$$

We define on $\mathbb N*\mathbb N $ an antimultiplicative involution by $\bar{e}=e$, $\bar{\alpha}=\beta$ and $\bar{\beta}=\alpha$. This involution extends by antilinearity into an involution of $A$, still denoted $*$.

Let $E_n\subset A$ be the linear space generated by the elements of $\mathbb N*\mathbb N$ having length $\leq n$.

Let $l^2(\mathbb N*\mathbb N)$ be the completion of $A$ with respect to the 2-norm. Observe that the elements of $\mathbb N*\mathbb N$ (viewed as elements of $A$) form an orthonormal basis of $l^2(\mathbb N*\mathbb N)$ .

Let $\tau_0:x\to <x(e),e>$ be the canonical state on $B(l^2(\mathbb N*\mathbb N))$. We define $S,T\in B(l^2(\mathbb N*\mathbb N))$ by linearity and by $S(x)=\alpha x,T(x)=\beta x$ for any $x\in\mathbb N*\mathbb N$.

\medskip

{\bf Reminder.} To any Hilbert space $H$ we can associate (cf. \cite{vdn}) the full Fock space $F(H)$: if $\{ h_i\}_{i\in I}$ is an orthonormal basis of $H$, then $\{h_{i_1}\otimes h_{i_2}\otimes\ldots\otimes h_{i_k}|k\geq 0\}$ is an orthonormal basis of $F(H)$ (at $k=0$, the corresponding vector is the vacuum one).

Let $\{f_i\}_{i\in I}$ be the generators of the free monoid $\mathbb N^{*I}$. Then the canonical basis of $l^2(\mathbb N^{*I})$ is the family $\{\delta_m\}$, with $m\in\mathbb N^{*I}$, of the form $m=f_{i_1}\ldots f_{i_k}$, and so:
$$l^2(\mathbb N^{*I})\simeq F(H)\quad,\quad\delta_{f_{i_1}\ldots f_{i_k}}\to h_{i_1}\otimes h_{i_2}\otimes\ldots\otimes h_{i_k}$$

The creation operator $l(h_i)$ corresponds in this way to $\lambda_{N^{*I}}
(f_i)$, where $\lambda_{N^{*I}}$ is the left regular representation (by isometries!) of the monoid $\mathbb N^{*I}$.

For $I=\{1,2\}$ we have $S=\lambda_{N*N}(\alpha )$ and $T=\lambda_{N*N} (\beta )$, so by identifying $l^2(\mathbb N*\mathbb N)=F(H)$, with
$H$ having orthonormal basis $\{ h_1,h_2\}$, we have:
$$S=l(h_1)\quad,\quad T=l(h_2)$$

With these conventions, we have the following result:

\begin{lemma} 
We define $\odot:\mathbb N*\mathbb N\times\mathbb N*\mathbb N\to A$ by: 
$$x\odot y=\sum_{x=ag,y=\bar{g}b}ab$$
\begin{enumerate}
\item $\odot$ extends by linearity into an associative multiplication of $A$.

\item If $P:(A,+,\cdot )\to (B(l^2(\mathbb N*\mathbb N)),+,\circ )$ is the $*$-morphism given by $\alpha\to S+T^*$ and $J:A\to A $ is the map $f\to P(f)e$, the $(J-Id)E_n\subset E_{n-1}$ for any $n$.

\item $J$ is an isomorphism of $*$-algebras $(A,+,\cdot)\simeq(A,+,\odot)$.
\end{enumerate}
\end{lemma}

\begin{proof}
(1) Observe first that $\odot$ is well-defined, the sum being finite. Let us prove now that $\odot$ is associative. Let $x,y,z\in\mathbb N*\mathbb N$. Then:
$$(x\odot y)\odot z=\sum_{\{g,a,b\in\mathbb N*\mathbb N|x=a\bar{g},y=gb\}} ab\odot z=\sum_{\{g,h,a,b,c,d\in\mathbb N*\mathbb N|x=a\bar{g},y=gb,ab=ch,z=\bar{h}d\}}cd$$

Now observe that for $a,b,c,h\in\mathbb N*\mathbb N$ the equality $ab=ch$ is equivalent to $b=uh,c=au$ with $u\in\mathbb N*\mathbb N$, or to $a=cv,h=vb$ with $v\in\mathbb N*\mathbb N$. Thus, we have:
$$(x\odot y)\odot z=\sum_{\{g,h,a,d,u\in\mathbb N*\mathbb N,x=a\bar{g},y=guh,z=\bar{h}d\}}aud+\sum_{\{g,b,c,d,v\in\mathbb N*\mathbb N ,x=cv\bar{g},y=gb,z=\overline{b}\bar{v}d\}}cd$$

A similar computation shows that $x\odot (y\odot z)$ is given by the same formula.

(2) Let $f\in A$. One can easily show that $P(\alpha )f=(S+T^*)f=\alpha\odot f$. Thus $J(\alpha g)=P(\alpha )J(g)=\alpha\odot J(g)=J(\alpha )\odot J(g)$ for any $g\in A$, and the same argument gives $J(\beta g)=J(\beta )\odot J(g)$,
for any $g\in A$. Now $(A,+,\cdot )$ being generated by $\alpha$ and $\beta$, we conclude that $J$ is a morphism of algebras, as follows:
$$J:(A,+,\cdot )\to (A,+,\odot )$$

We prove now by recurrence on $n\geq1$ that $(J-Id)E_n\subset E_{n-1}$. At $n=1$ we have $J(\alpha )=\alpha$, $J(\beta)=\beta$ and $J(e)=e$, and since $E_1$ is generated by $e,\alpha ,\beta$, we have $J=Id$ on $E_1$. Now assume that this is true for $n$, and let $k\in E_{n+1}$. We write $k=\alpha f+\beta g+h$ with $f,g,h\in E_n$ (observe that this decomposition is not unique). Then:
\begin{eqnarray*}
(J-Id)k
&=&J(\alpha f+\beta g+h)-(\alpha f+\beta g+h)\\
&=&[(S+T^*)J(f)+(S^*+T)J(g)+J(h)]-[Sf+Tg+h]\\
&=&S(J(f)-f)+T(J(g)-g)+T^*J(f)+S^*J(g)+(J(h)-h)
\end{eqnarray*}

By using the recurrence assumption, applied to $f,g,h$ we find that $E_n$ contains all the terms of the above sum, and so contains $(J-Id)k$, and we are done.

(3) Here we have to prove that $J$ preserves the involution $*$, and that it is bijective. We have $J*=*J$ on the generators $\{ e,\alpha ,\beta\}$ of $A$, so $J$ preserves the involution. Also, by (2), the restriction of $J-Id$ to $E_n$ is nilpotent, so $J$ is bijective. 
\end{proof}

We can now start the representation theory work, as follows:

\begin{lemma}
Let $(G,u)$ be a compact matrix quantum group, and let $\Psi_G: (A,+,\odot )\to G$ be the unique morphism defined by $\alpha\to\chi (u)$, $\beta\to\chi(\bar{u})$, cf. Lemma 3 (3).

Let $n\geq 1$ and assume that $\Psi_G(x)$ is the character of an irreducible representation $r_x$ of $G$, for any $x\in\mathbb N*\mathbb N$ having length $\leq n$. Then $\Psi_G(x)$ is the character of a (non-null) representation of $G$, for any $x\in\mathbb N*\mathbb N$ of length $n+1$.
\end{lemma}

\begin{proof}
At $n=1$ this is clear. Assume $n\geq2$, and let $x\in\mathbb N*\mathbb N$ of length $n+1$. If $x$ contains a $\geq 2$ power of $\alpha$ or of $\beta$, for instance if $x=z\alpha^2y$, then we can set $r_x=r_{z\alpha}\otimes r_{\alpha y}$ and we are done. So, assume that $x$ is an alternating product of $\alpha$ and $\beta$. We can assume that $x$ begins with $\alpha$. Then $x=\alpha\beta\alpha y$, with $y\in\mathbb N*\mathbb N$ being of length $n-2$. 

Observe that the equality $\Psi_G(\bar{z})=\Psi_G(z)^*$ holds on the generators $\{ e,\alpha ,\beta\}$ of $\mathbb N*\mathbb N$, so it holds for any $z\in\mathbb N*\mathbb N$. Thus if $<,>$ is the scalar product on $G$ associated to the Haar measure (see the Reminders), we have: 
\begin{eqnarray*}
<\chi (r_\alpha\otimes r_{\beta\alpha y}),\chi (r_{\alpha y})>
&=&<\chi (r_{\beta\alpha y}),\chi (r_\beta\otimes r_{\alpha y})>\\
&=&<\chi (r_{\beta\alpha y}),\Psi_G(\beta\odot \alpha y )>\\
&=&<\chi (r_{\beta\alpha y}),\Psi_G(\beta\alpha y) +\Psi_G(y)>\\
&=&<\chi (r_{\beta\alpha y}),\chi (r_{\beta\alpha y})+\chi (r_y)>\\
&\geq&1
\end{eqnarray*}

Since $r_{\alpha y}$ is irreducible, we have $r_{\alpha y}\subset r_\alpha\otimes r_{\beta\alpha y}$. Thus $\chi (r_\alpha\otimes r_{\beta\alpha y})-\chi (r_{\alpha y})=\Psi_G(\alpha\odot \beta\alpha y
-\alpha y)=\Psi_G(x)$ is the character of a representation of $G$.
\end{proof}

Let $(G,u)$ be a compact matrix quantum group, and set $f_x=\Psi_G(x)$ for any $x\in\mathbb N*\mathbb N$, with notations from Lemma 4. Then the family $\{f_x|x\in\mathbb N*\mathbb N\}$ satisfies $f_e=1$, $f_\alpha=\chi(u)$, $f_\beta=\chi(\bar{u})$ and $f_xf_y=\sum_{x=ag,y=\bar{g}b}f_{ab}$. We want to prove:
\begin{enumerate}
\item Theorem 1 (1): for $G=A_u(F)$, these elements $f_x$ are precisely the characters of the irreducible representations of $A_u(F)$.

\item Theorem 2: if these elements $f_x$ are the characters of the irreducible representations of $G$, then we have $G_{full}\sim_{sim}A_u(F)$, for a certain matrix $F$.
\end{enumerate}

It is convenient to consider, for any $n\in\mathbb N$, any $F\in GL(n,\mathbb C)$, and any compact matrix quantum group $(G,u)$ with $u,F\bar{u}F^{-1}$ unitaries, the following diagram:
$$\begin{matrix}
\mathbb N*\mathbb N &\ &\  \mathbb N *\mathbb N  &\ &\ A_u(F)\\
&\ \\
\cap &\ &\ \cap &\ \nearrow\Psi_u &\ \downarrow\Phi\\
&\ \\
C<X,X^*>=(A,+,\cdot )&\ \displaystyle{\mathop{\longrightarrow}^J}&\ (A,+,\odot )&\
\displaystyle{\mathop{\longrightarrow}^{\Psi_G}}&G_{full}\cr
&\ \\
P\downarrow &\ &\ \tau\downarrow &\ (\star )&\ \downarrow h\cr
&\ \\
B(l^2(\mathbb N *\mathbb N))&\ \displaystyle{\mathop{\longrightarrow}^{\tau_0}}&\ \mathbb C&\
\displaystyle{\mathop{\longrightarrow}^{Id_\mathbb C }}&\mathbb C
\end{matrix}$$

To be more precise, the various objects and arrows are as follows:

-- $\tau_0,J,P$ are as above, and $\tau (f)=f(e)=$ coefficient of $e$ inside $f$. 

-- $G_{full}$ is the full version of $G$ (see Reminders).

-- $\Phi$ is the canonical surjection, defined by the universality property of $A_u(F)$.

-- $\Psi_G$ (resp. $\Psi_u$) is the unique $*$-morphism (see Lemma 3 (3)) which maps $\alpha$ on the character of the fundamental representation of $G_{full}$ (resp. of $A_u(F)$).

-- $h$ is the Haar measure of $G$. 

Note that the square on the left and the triangle clearly commute. Observe also that the inclusions $\mathbb N*\mathbb N\subset A$ (see Notations above) do not commute with $J$. We have:

\begin{proposition}
Let $(G,u)$ be a compact matrix quantum group, with $u$ and $F\bar{u}F^{-1}$ assumed to be unitaries. The following are equivalent:
\begin{enumerate}
\item The irreducible representations of $G$ can be indexed by $\mathbb N*\mathbb N$, with $r_e=1,r_{\alpha}=u,r_{\beta}=\bar{u}$ and $r_x\otimes r_y=\sum_{x=ag,y=\bar{g}b}r_{ab}$.

\item The square $(\star )$ in the above diagram commutes.

\item $\chi (u)/2$ is a circular variable in $(G,h)$.

\item The $*$-moments of $\chi (u)/2\in (G,h)$ are smaller than the respective $*$-moments of a circular variable, i.e. $\mu_{\chi (u)/2}(M)\leq \mu_c(M)$ for any noncommutative monomial $M$ (where $\mu_c$ is the $*$-distribution of the circular variable).
\end{enumerate}
In addition, if these conditions are satisfied, then $\Phi :A_u(F)\to G_{full}$ is a $*$-isomorphism.
\end{proposition}

{\bf Note.} We will see in the next section that the last condition is equivalent to (1-4).

\begin{proof}
$(1)\implies(2)$ It is clear that $\mathbb N*\mathbb N$ is an orthonormal system in $((A,+,\odot ),\tau )$. Thus $\Psi_G(\mathbb N*\mathbb N)=\{\chi (r_x)|x\in\mathbb N*\mathbb N\}$ is an orthonormal system in $(G_s,h)$, so ($\star$) commutes.

$(2)\implies(3)$ Let us first remark that if ($\star$) commutes then $h\Psi_GJ=\tau_0P$. Now by identifying $(\mathbb C<X,X^*>,+,\cdot )=(A,+,\cdot )$ we have (see Reminders (G)):

-- the $*$-distribution of $\chi (u)\in (G,h)$ is the functional $h\Psi_GJ$.

-- the $*$-distribution of $S+T^*\in(B(l^2(\mathbb N*\mathbb N)),\tau_0)$ is the functional $\tau_0P$.

On the other hand, the identifications made in the beginning of this section show that $(S+T^*)/2$ has the same $*$-distribution as the variable $(l(h_1)+l(h_2)^*)/2$ on the full Fock space, which is the standard example of circular variable (see the $1^{st}$ section of \cite{vdn}).

$(3)\implies(4)$ is trivial. Let us prove now $(4)\implies(1)$. By identifying as usual $(\mathbb C<X,X^*>,+,\cdot )=(A,+,\cdot )$, the noncommutative monomials in $X,X^*$ correspond to the elements of $\mathbb N*\mathbb N\subset A$. Thus, the assumption on the $*$-moments of $\chi (u)/2$ reads:
$$(i)\,\,\,\,\,\,\,\,\,\,\,\,\, h\Psi_GJ\leq \tau J\,\,\,\, {\rm on}\,\,\,\,  \mathbb N*\mathbb N$$

We prove by recurrence on $n\geq 0$ that for any $z\in \mathbb N*\mathbb N$ having length $n$, $\Psi_G(z)$ is the character of an irreducible representation $r_z$ of $G$. At $n=0$ we have $\Psi_G(e)=1$. So, assume that our claim holds at $n\geq 0$, and let $x\in\mathbb N*\mathbb N$ having length $n+1$.

Lemma 3 (2) gives $J(x)=x+z$ with $z\in E_n$. Let $A_N\subset A$ be the set of functions $f$ such that $f(x)\in\mathbb N$ for any $x\in\mathbb N*\mathbb N$. Then $J(\alpha ),J(\beta )\in A_N$, so by multiplicativity $J(\mathbb N*\mathbb N)\subset A_N$. In particular, $J(x)\in A_N$. Thus there exist $m(z)\in\mathbb N$ such that:
$$J(x)=x+\sum_{z\in\mathbb N*\mathbb N,l(z)\leq n}m(z)z$$

Let us compute, by using this formula, $h\Psi_GJ(x\bar{x})$ and $\tau J(x\bar{x})$: 

(a) It is clear that for $a,b\in\mathbb N*\mathbb N$ we have $\tau(a\odot\bar{b})=\delta_{a,b}$. Thus:
$$\tau J(x\bar{x})=\tau((x+\sum m(z)z)\odot (\bar{x}+\sum m(z)\bar{z}))=1+\sum m(z)^2$$

(b) By recurrence and by Lemma 4, $\Psi_G(x)$ is the character of a representation $r_x$ of $G$. Thus  $\Psi_GJ(x)$ is the character of $r_x+\sum_{z\in\mathbb N*\mathbb N,l(z)\leq n}m(z)r_z$. Now by using the orthogonality formulae for the characters, we obtain:
$$h\Psi_GJ(x\bar{x})\geq h(\chi (r_x)\chi (r_x)^*)+\sum m(z)^2$$

By $(i)$, (a), (b) we conclude that $r_x$ is irreducible, which finishes the recurrence. Finally, the fact that the $r_x$ are distinct comes from $(i)$. Indeed, $\mathbb N*\mathbb N$ being an orthonormal basis of $((A,+,\odot ),\tau )$, for any $x,y\in\mathbb N*\mathbb N$, $x\neq y$ we have $\tau (x\odot\bar{y})=0$, and so: 
$$h(\chi (r_x\otimes\bar{r_y}))=h\Psi_GJ(x\bar{y})\leq \tau
J(x\bar{y})=\tau (x\odot\bar{y})=0$$

It remains to prove that (1) implies the last condition. We prove by recurrence on $n\geq 0$ that for any $x\in\mathbb N*\mathbb N$ having length $n$, $\Psi_u(x)$ is the character of an irreducible representation $p_x$ of $A_u(F)$. At $n=0$ this is trivial. So, assume that this holds for $n$, and let $x\in\mathbb N*\mathbb N$ having length $n+1$. By Lemma 4, $\Psi_u(x)$ is the character of a representation $p_x$. Since $\Phi$ maps $p_x\to r_x$, which is irreducible, $p_x$ is irreducible too, and we are done. 

Summarizing, our surjection $\Phi$ maps the classes of irreducible representations of $A_o(F)$ to the classes of  irreducible representations of $G_{full}$. We can conclude by using a standard argument (Theorem 5.7 in \cite{wo2}). Indeed, if $\{ c_i\}$ is a basis of $A_u(F)_s$ formed by coefficients of irreducible representations, then $\{ \Phi (c_i)\}$ is a basis of $G_s$ formed by coefficients of irreducible representations. Thus, $\Phi :A_u(F)\to G_{full}$ is bijective.
\end{proof}

\section{Representations of $A_u(F)$}

In this section we prove Theorem 1 and Theorem 2. We will see that  Proposition 3 implies easily Theorem 2, as well as (modulo a free probability result) Theorem 1 for the matrices satisfying $F\bar{F}\in\mathbb RI_n$. Theorem 1 will be afterwards proved for the arbitrary matrices $F\in GL(n,\mathbb C)$, by using the results obtained for $F=I_n$.

\begin{proof}(of Theorem 2) Let $(G,u)$ be a compact matrix quantum group, whose irreducible representations are labelled by $\mathbb N*\mathbb N$, such that $r_e=1$, $r_\alpha =u$, $r_\beta =\bar{u}$ and:
$$r_x\otimes r_y=\sum_{x=ag,y=\bar{g}b}r_{ab}$$ 

We can assume (modulo similarity) that $u\in M_n(G)$ is unitary. Thus, there exists $F\in GL(n,\mathbb C)$ such that $F\bar{u}F^{-1}$ is unitary, and the result follows from Proposition 3. 
\end{proof}

\begin{proof}(of Theorem 1, case $F\bar{F}\in\mathbb RI_n$) Let $u,v$ be the fundamental representations of $A_u(F),A_o(F)$, and let $z\in C(\mathbb T)$ be the map $x\to x$. Let $G$ be the $C^*$-subalgebra of $C(\mathbb T)*_{red}A_o(F)$ generated by the coefficients of $zv=(zv_{ij})_{ij}$. We have then:

-- $\chi (v)/2$ is semicircular inside $A_o(F)$ (cf. Proposition 1).

-- $z$ is a Haar unitary in $C(\mathbb T)$ (clear).

-- $\chi (v)/2$ and $z$ are $*$-free in $C(\mathbb T)*_{red}A_o(F)$ (cf. Reminders (F)).

These three conditions show that the product $z\chi (v)/2$ is circular in $C(\mathbb T)*_{red}A_o(F)$ (this is a well-known version of Voiculescu's theorem on the polar decomposition of the circular variables \cite{vdn}, see for instance \cite{ba1} or \cite{nsp}). But $z\chi (v)=\chi (zv)$ is the character of the fundamental representation of  $(G,zv)$, and so we can apply Proposition 3:

-- The fact that the condition (3) there implies the last condition there shows that we have $A_u(F)=G_{full}$, and so $A_u(F)_{red}=G$, which gives Theorem 1 (4).

-- The fact that we have $(3)\implies(1)$ there classifies the representations of $G_{full}=A_u(F)$, and we obtain in this way Theorem 1 (1,2,3).
\end{proof}

{\bf Remark.} We could have proved Theorem 1 in the case $F\bar{F}\in\mathbb RI_n$ in the following way. Consider the compact quantum group $G\subset C(\mathbb T)*_{red}A_o(F)$ generated by the coefficients of $zv$, with $v$ being the fundamental representation of $A_o(F)$. By combining the theory of representations of $A_o(F)$ from \cite{ba2} and the theory of representations of the free products from \cite{wa1}, we can classify the representations of $G$. Then, we can apply Theorem 2, in order to conclude that we have $A_u(F)_{red}=G$. Note however that this proof does not provide us with any tool for investigating the general case, $F\in GL(n,\mathbb C)$ arbitrary.

\medskip

We discuss in what follows the proof of Theorem 1 in the general case, where $F\in GL(n,\mathbb C)$ is arbitrary. Let $u$ be the fundamental representation of $A_u(F)$. We must compute the $*$-moments of the character $\chi (u)$, i.e. the following numbers:
$$h(\chi (u)^{a_1}\chi (u)^{*b_1}\chi (u)^{a_2}\ldots)=dim(Mor(1,u^{\otimes a_1}\otimes \bar{u}^{\otimes b_1}\otimes u^{\otimes
a_2}\otimes\ldots))$$ 

By using the fact that $\mathbb N*\mathbb N$ is a free monoid:

-- we associate to any Hilbert space $H$ a family of Hilbert spaces $\{
H_x\}_{x\in\mathbb N*\mathbb N}$ in the following way: $H_e=\mathbb C$, $H_\alpha =H$, $H_\beta=\bar{H}$, and
$H_{ab}=H_a\otimes H_b$, $\forall\, a,b\in\mathbb N*\mathbb N$. 

-- we define a family $\{u^x\}_{x\in\mathbb N*\mathbb N}$ of unitary representations of $A_u(F)$ in the following way: $u^e=1$, $u^\alpha=u$, $u^\beta=F\bar{u}F^{-1}$ (acting on $\bar{\mathbb C}^n$) and $u^{ab}=u^a\otimes u^b$, $\forall\,a,b\in\mathbb N*\mathbb N$. 

Observe that we have $u^x\in B(\mathbb C^n_x)\otimes A_u(F)$, for any $x$. With these conventions, the $*$-moments of $\chi (u)$ are the following numbers:
$$\{\dim(Mor(1,u^k))|k\in \mathbb N*\mathbb N\}$$

We will compute these numbers by using Theorem 1.3. from \cite{wo3}, as follows:

\begin{lemma}
Let $F\in GL(n,\mathbb C)$, set $H=\mathbb C^n$, and consider the linear maps $E_1:H_e\to H_{\alpha\beta}$ given by $1\to\sum F(e_i)\otimes \bar{e}_i$ and $E_2:H_e\to
H_{\beta\alpha}$ given by $1\to\sum\bar{e}_i\otimes\bar{F}^{-1}(e_i)$. 
\begin{enumerate}
\item $E_1\in Mor(1,u^{\alpha\beta})$ and $E_2\in Mor(1,u^{\beta\alpha})$.

\item $(E_2^*\otimes Id_{H_\beta})(Id_{H_\beta}\otimes E_1)\in \mathbb C Id_{H_\beta}$ and $(E_1^*\otimes Id_{H_\alpha})(Id_{H_\alpha}\otimes E_2)\in\mathbb C Id_{H_\alpha}$.

\item For $r,s\in\mathbb N*\mathbb N$ let $Mor(r,s)\subset B(H_r,H_s)$ be the sets of linear combinations of (composable) products of linear maps of the form $Id_{H_k}$ or $Id_{H_k}\otimes E_1\otimes Id_{H_p}$ or $Id_{H_k}\otimes E_2\otimes Id_{H_p}$ or $Id_{H_k}\otimes E_1^*\otimes
Id_{H_p}$ or $Id_{H_k}\otimes E_2^*\otimes Id_{H_p}$. Then the concrete monoidal $W^*$-category of representations of $A_u(F)$ is the completion, in the sense of \cite{wo3}, or the following concrete monoidal $W^*$-category: 
$$Z(F)=\{\mathbb N*\mathbb N,\cdot ,\{ H_r\}_{r\in\mathbb N*\mathbb N},\{ Mor(r,s)\}_{r,s\in\mathbb N*\mathbb N}\}$$
\end{enumerate}
\end{lemma}

\begin{proof}
(1) If $w\in M_n(B)$ is a unitary over a $*$-algebra and $\zeta =\sum e_i\otimes e_i\in \mathbb C^n\otimes \mathbb C^n$ then:
$$(w_{13}\bar{w}_{23})(\zeta\otimes 1)=\sum e_i\otimes e_k\otimes w_{ia}w_{ka}^*=\sum e_i\otimes e_k\otimes \delta_{ik}1=(\zeta\otimes 1)$$

By using this formula, we obtain:

-- with $B=A_u(F)$ and $w=u$ this shows that $(1\otimes F)\zeta =\sum e_i\otimes Fe_i$ is a fixed vector of $(1\otimes F)(u\otimes\bar{u})(1\otimes F)^{-1}=u\otimes (F\bar{u}F^{-1})$.

-- with $B=A_u(F)$ and $w=F\bar{u}F^{-1}$ this shows that $(1\otimes\bar{F}^{-1})\zeta =\sum e_i\otimes\bar{F}^{-1}e_i$ is a fixed vector of $(1\otimes\bar{F})^{-1}(w\otimes\bar{w})(1\otimes\bar{F})=(F\bar{u}F^{-1})\otimes u$.

On the other hand, by definition of $u^x$ we have $u^{\alpha\beta}=(id_{\mathbb C^n}\otimes \phi\otimes id_{A_u(F)} )(u\otimes (F\bar{u}F^{-1}))$, and we have as well $u^{\beta\alpha}=(\phi\otimes id_{\mathbb C^n}\otimes id_{A_u(F)})((F\bar{u}F^{-1})\otimes u)$, where $\phi :B(\mathbb C^n)\to B(\bar{\mathbb C}^n)$ is the canonical isomorphism. Thus, we obtain (1).

(2) is an easy computation. Regarding now (3), observe that $Z(F)$ is by definition a concrete monoidal $W^*$-category. Let $j:H_\alpha\to H_\beta$ be the antilinear map defined by $e_i\to F(\bar{e}_i)$. Then (with the notations from \cite{wo3}, page 39) we have $t_j=E_1$ and $(\bar{t}_j)^*=t_{j^{-1}}=E_2$, so $t_j\in Mor(e,\alpha\beta)$ and $\bar{t}_j\in Mor(1,\beta\alpha)$, so $\bar{\alpha}=\beta$ in $Z(F)$. By \cite{wo3}, the universal $Z(F)$-admissible pair is a certain compact quantum group $(G,v)$. 

Since by (1), $(A_u(F),u)$ is a $Z(F)$-admissible pair, we have a $C^*$-morphism $f:G\to A_u(F)$ such that $(id\otimes f)(v)=u$. On the other hand, the universal property of $A_u(F)$ gives a $C^*$-morphism $g:A_u(F)\to G$ such that $(id\otimes g)(u)=v$. Thus $(G,v)=(A_u(F),u)$.
\end{proof}

The proof of the following lemma is similar to the proof of Lemma 2:

\begin{lemma}
Let $I(p)=Id_{H_p}$ and $V_i(p,q)=I(p)\otimes E_i\otimes I(q)$, for $i=1,2$. 
\begin{enumerate}
\item Any morphism of $Z(F)$ is a linear combination of maps of the form $I(.)$ or of the form $V_\cdot(.,.)\circ\ldots\circ V_\cdot (.,.)\circ V_\cdot (.,.)^*\circ\ldots\circ V_\cdot (.,.)^*$.

\item The maps of the form $M\otimes I(\alpha )\otimes N$ with $M\in Mor(e,x)$, $N\in Mor(e,y)$ and $x\alpha y=k$ generate $Mor(\alpha ,k)$. lso, the maps of the form $M\otimes I(\beta )\otimes N$ with $M\in Mor(e,x)$, $N\in Mor(e,y)$ and $x\beta y=k$ generate $Mor(\beta ,k)$.

\item The maps of the form $(I(\alpha )\otimes M\otimes I(\beta )\otimes N)\circ E_1$ with $M\in Mor(e,x)$, $N\in Mor(e,y)$ and $\alpha x\beta y=k$, or of the form $(I(\beta )\otimes M\otimes I(\alpha )\otimes N)\circ E_2$ with $M\in Mor(e,x)$, $N\in Mor(e,y)$ and $\beta x\alpha y=k$ generate $Mor(e,k)$.
\end{enumerate}
\end{lemma}

As a final ingredient now, we have:

\begin{proposition}
For any $k\in\mathbb N*\mathbb N$ we define $Z_k(F)\subset Mor(e,k)$ by $Z_e(F)=1$, $Z_\alpha(F)=Z_\beta (F)=\emptyset$ and (by recurrence)
$$Z_k(F)=\cup_{k=\alpha x\beta y}\{(I(\alpha )\otimes M\otimes I(\beta )\otimes N)\circ E_1|M\in Z_x(F), N\in Z_y(F)\}$$
if $k$ begins with $\alpha$ and 
$$Z_k(F)=\cup_{k=\beta x\alpha y}\{ (I(\beta )\otimes M\otimes
I(\alpha )\otimes N)\circ E_2\mid M\in Z_x(F), N\in Z_y(F)\}$$
if $k$ begins with $\beta$. Then:
\begin{enumerate}
\item $Z_k(F)$ generates $Mor(e,k)$.

\item $Z_k(I_n)$ is a basis of $Mor(e,k)$ (for $F=I_n$). 
\end{enumerate}
\end{proposition}

\begin{proof}
Here (1) follows from Lemma 6 (3). In order to prove now (2), we recall that $H_k$ is a certain tensor product between $H$ and its conjugate $\bar{H}$. Let $\psi:H\to\bar{H}$ be the isometry given by $e_i\to\bar{e}_i$. By identifying $\bar{H}$ with $H$ via $\psi$, we obtain an isometry $\psi_k:H_k\to H^{\otimes l(k)}$, where $l(k)$ is the length of the word $k$.

It is clear that $\psi_k$ maps $\{ X(1)|X\in Z_k(I_n)\}$ to a subset of $\{ X(1)|X\in W_{l(k)}(I_n)\}$. By using Proposition 2, the set $\{ X(1)|X\in W_{l(k)}(I_n)\}$ is formed by linearly independent vectors in $H^{\otimes l(k)}$. Thus $\{ X(1)|X\in Z_k(I_n)\}$ is formed by linearly independent vectors, and so $Z_k(I_n)$ is a basis of $Mor(e,k)$.
\end{proof}

We can now finish the proof of Theorem 1:

\begin{proof}(of Theorem 1, general case) We recall that $F\in GL(n,\mathbb C)$ is now arbitrary. By construction of $Z_k(F)$ we have $Card(Z_k(F))=C_k$, where $\{ C_x\}_{x\in\mathbb N*\mathbb N}$ are the numbers defined by $C_e=1$, $C_\alpha=C_\beta=0$ and $C_k=\sum_{k=\alpha x\beta y}C_xC_y+\sum_{k=\beta x\alpha y}C_xC_y$ for $k\in\mathbb N*\mathbb N$.

Let $u(F)$ be the fundamental representation of $A_u(F)$. By using Lemma 5 (3) and Proposition 4 (1), for any $k\in\mathbb N*\mathbb N$ we have:
$$\dim(Mor(1,u(F)^k))=\dim(Mor(e,k))\leq Card(Z_k(F))=C_k$$

Proposition 4 (2) tells us that for $F=I_n$ we have equality. Thus:
$$\dim(Mor(1,u(F)^k))\leq\dim(Mor(1,u(I_n)^k))$$

But the numbers $\dim(Mor(1,u(F)^k))$ are the $*$-moments of the character of $u(F)$, and the numbers $\dim(Mor(1,u(I_n)^k))$ are the $*$-moments of the character of $u(I_n)$, and so the $*$-moments of a circular variable (by Theorem 1 (3) applied to the case $F=I_n$, that we have already solved). We can conclude by using Proposition 3.
\end{proof}

{\bf Remark.} It is clear now that for any matrix $F$, the set $\{ X(1)|X\in Z_k(F)\}$ is a basis of the space of fixed vectors for the representation 
$u(F)^k$ of $A_u(F)$. If $F\bar{F}\in\mathbb R I_n$, the map $\psi_k$ constructed in the proof of Proposition 4 allows us to identify this basis with a part of the basis of fixed vectors for the representation $u^{\otimes l(k)}$ of $A_u(F)$, constructed at the end of section 2 (in fact, this identification is the one coming from $A_u(F)\to A_o(F)$).

It is easy to determine the noncrossing pairings of $\{1,\ldots,l(k)\}$ which correspond to the fixed vectors of the representation $u^k$ of $A_u(F)$. In fact, we can obtain such a basis by writing $k=x_1\ldots x_{l(k)}$ with $x_i\in\{\alpha ,\beta\}$ and by associating to any noncrossing pairing $P=P_1\coprod\ldots\coprod P_{l(k)/2}$ of $\{1,\ldots,k\}$ with parts of the form $P_s=\{ i_s,j_s\}$ such that $(x_{i_s},x_{j_s})$ is equal to $(\alpha ,\beta )$ or to $(\beta ,\alpha )$ for any $s$ the following vector:
$$v=\sum_{1\leq s_1,\ldots,s_{l(k)}\leq
l(k)}F_{s_{j_1}s_{i_1}}\ldots F_{s_{j_{l(k)/2}}s_{i_{l(k)/2}}}e_{s_1}\otimes\ldots\otimes e_{s_{l(k)}}$$

\section{The case $n=2$}

We recall that given $\mu\in [-1,1]-\{ 0\}$, the compact matrix quantum group $S_\mu U(2)$ is constructed with generators $\alpha ,\gamma$ and the following relations:
$$\alpha^*\alpha +\gamma^*\gamma =1,\,\alpha\alpha^*
+\mu^2\gamma\gamma^*=1,\,\gamma\gamma^*=\gamma^*\gamma ,\, \mu\gamma\alpha =\alpha\gamma
,\,\mu\gamma^*\alpha =\alpha\gamma^*$$

We have $S_1U(2)=C(SU(2))$ (see \cite{wo1} and the formulae (1.33) in  \cite{wo2}). We have:

\begin{proposition}
$A_o\begin{pmatrix}0&1\\-\mu^{-1}&0\end{pmatrix}=S_{\mu}U(2)$.
\end{proposition}

\begin{proof}
Let $\begin{pmatrix}u_{11}&u_{12}\\u_{21}&u_{22}\end{pmatrix}$ be the fundamental representation of $A_o\begin{pmatrix}0&1\\-\mu^{-1}&0\end{pmatrix}$. The relations defining $A_o\begin{pmatrix}0&1\\-\mu^{-1}&0\end{pmatrix}$ are $uu^*=u^*u=1$ and:
$$\begin{pmatrix}u_{11}&u_{12}\\u_{21}&u_{22}\end{pmatrix}=\begin{pmatrix}0&1\\-\mu^{-1}&0\end{pmatrix}
\begin{pmatrix}u_{11}^*&u_{12}^*\\u_{21}^*&u_{22}^*\end{pmatrix}
\begin{pmatrix}0&-\mu\\1&0\end{pmatrix}$$

By multiplying, $\begin{pmatrix}u_{11}&u_{12}\\u_{21}&u_{22}\end{pmatrix}=\begin{pmatrix}u_{22}^*&-\mu
u_{21}^*\\-\mu^{-1}u_{12}^*&u_{11}^*\end{pmatrix}$. With $\alpha=u_{11}$ and $\gamma=u_{21}$ we have $u=\begin{pmatrix}\alpha&-\mu\gamma^*\\\gamma&\alpha^*\end{pmatrix}$ and the relations defining $A_o\begin{pmatrix}0&1\\-\mu^{-1}&0\end{pmatrix}$ are $uu^*=u^*u=1$, i.e.:
$$\begin{pmatrix}\alpha&-\mu\gamma^*\\\gamma&\alpha^*\end{pmatrix}
\begin{pmatrix}\alpha^*&\gamma^*\\-\mu\gamma&\alpha\end{pmatrix}
=\begin{pmatrix}\alpha^*&\gamma^*\\-\mu\gamma&\alpha\end{pmatrix}
\begin{pmatrix}\alpha&-\mu\gamma^*\\\gamma&\alpha^*\end{pmatrix}
=\begin{pmatrix}1&0\\0&1\end{pmatrix}$$

By computing we obtain the relations defining $S_\mu U(2)$.
\end{proof}

We will need the following result, regarding the general case, $n\in\mathbb N$:

\begin{proposition}
For any $n\in\mathbb N$, $F\in\ GL(n,\mathbb C),\,\lambda\in\mathbb C^*$ and $V,W\in U(n)$ we have:
\begin{enumerate}
\item  $A_o(F)\sim_{sim} A_o(\lambda VFV^t)$ (if $F$ satisfies $F\overline{F}\in\mathbb R I_n$).

\item $A_u(F)\sim_{sim} A_u(\lambda VFW)$.
\end{enumerate}
\end{proposition}

\begin{proof}
(1) Let $u,v$ be the fundamental representations of $A_o(F)$ and of  $A_o(\lambda VFV^t)$. Then $v=(\lambda VFV^t)\bar{v}(\lambda
VFV^t)^{-1}$ is unitary $\implies$ $v=VFV^t\bar{v}\bar{V}F^{-1}V^*$ is unitary. Thus $V^*vV=F\overline{V^*vV}F^{-1}$ is unitary, so we can define $f:A_o(F)\to A_o(\lambda VFV^t)$ by the universality property and by $(Id\otimes f)(u)=V^*vV$. The same arguments show that we have a morphism 
$g: A_o(\lambda VFV^t)\to A_o(F)$ with $(Id\otimes g)(v)=VuV^*$, so $f,g$ are inverse bijections.

(2) Let $u,v$ be the fundamental representations of $A_u(F)$ and of  $A_u(\lambda VFW)$. Then $v$ and $(\lambda VFW)\bar{v}(\lambda VFW)^{-1}$ are unitaries, so $\bar{W}vW^t$ and $FW\bar{v}W^*F^{-1}$ are unitaries, so $\bar{W}vW^t$ and $F\overline{\bar{W}vW^t}F^{-1}$ are unitaries. Thus we can define $f:A_u(F)\to A_u(\lambda VFW)$ by the universal property and by $(Id\otimes f)(u)=\overline{W}vW^t$. Similarly, we have a morphism $g: A_u(\lambda VFW)\to A_u(F)$ with $(Id\otimes g)(v)=W^tu\bar{W}$, so $f,g$ are inverse bijections.
\end{proof}

Now back to the case $n=2$, we have the following result:

\begin{proposition}
For $\mu\in [-1,1]-\{0\}$ let $G_\mu\subset C(\mathbb T)*_{red}S_\mu U(2)$ be the $C^*$-algebra generated by the coefficients of $zu_\mu$, with $u_\mu$ being the fundamental representation of $S_\mu U(2)$. Then, modulo similarity, we have:
\begin{enumerate}
\item $\{ A_o(F)|F\in GL(2,\mathbb C)\}=\{ S_\mu U(2)|\mu\in[-1,1]-\{ 0\}\}$.

\item $\{ A_u(F)_{red}|F\in GL(2,\mathbb C )\} = \{ G_\mu|\mu\in [-1,1]-\{ 0\}\}$.
\end{enumerate}
\end{proposition}

\begin{proof}
(1) We must show that for $F\in GL(2,\mathbb C)$ such that $F\bar{F}\in\mathbb R$, there exists $V\in GL(2,\mathbb C)$, scalar multiple of a unitary, such that $\lambda VFV^t=\begin{pmatrix}0&1\\-\mu^{-1}&0\end{pmatrix}$ for a certain $\mu$ (cf. Propositions 5 and 6). So, let $F=\begin{pmatrix}x&y\\ z&t\end{pmatrix}$. If $x\neq 0$, let $\alpha$ be a solution of $\alpha^2x+\alpha (y+z)+t=0$ and $V=\begin{pmatrix}\alpha&1\\-1&\bar{\alpha}\end{pmatrix}$, which is a multiple of a unitary. Then $VFV^t$ has as first entry $\alpha^2x+\alpha (y+z)+t=0$, so we can assume $x=0$. In the case $F=\begin{pmatrix}0&y\\ z&t\end{pmatrix}$ we have $F\bar{F}=\begin{pmatrix}y\bar{z}&y\bar{t}\\ t\bar{z}&z\bar{y}+t\bar{t}\end{pmatrix}$. Since $F\overline{F}\in\mathbb R$, $t=0$ and $y\bar{z}=z\bar{y}$, so $F=y\begin{pmatrix}0&1\\k&0\end{pmatrix}$ with $k=z/y\in\mathbb R ^*$. If $|k|\geq 1$ we are done, and otherwise, we can take $V=\begin{pmatrix}0&-1\\1&0\end{pmatrix}$.

(2) It is enough to show that for any $F\in GL(2,\mathbb C)$ there exist $V,W\in U(2)$ such that $VFW\overline{VFW}\in\mathbb R$ (cf. (1) and Proposition 6). By polar decomposition we can assume $F> 0$, and by diagonalizing, we can further assume $F=\begin{pmatrix}x&0\\ 0&y\end{pmatrix}$ with $x,y>0$. But in this case we can take $V=\begin{pmatrix}0&-1\\1&0\end{pmatrix}$ and $W=1$. 
\end{proof}

We can further build on the above observations, at $F=I_2$, as follows:

\begin{lemma}
Let $\begin{pmatrix}a&b\\-\bar{b}&\bar{a}\end{pmatrix}$, $(z)$ and $(u)$ be the fundamental representations of $C(SU(2))$, $C(\mathbb T)$ and $A_u(I_2)$. We have then en embedding, as follows:
$$A_u(I_2)_{red}\subset C(\mathbb T)*_{red}C(SU(2))\quad,\quad u\to\begin{pmatrix}za&zb\\-z\bar{b}&z\bar{a}\end{pmatrix}$$
\end{lemma}

\begin{proof}
We have indeed $A_o\begin{pmatrix}0&1\\-1&0\end{pmatrix}=S_1U(2)=C(SU(2))$ (cf. Proposition 5), and also $A_u(I_2)=A_u\begin{pmatrix}0&1\\-1&0\end{pmatrix}$ (since $A_u(F)=A_u(\sqrt{F^*F})$ for any $F$), and by Theorem 1 we obtain an embedding $A_u\begin{pmatrix}0&1\\-1&0\end{pmatrix}\subset C(\mathbb T)*_{red}A_o\begin{pmatrix}0&1\\-1&0\end{pmatrix}$.
\end{proof}

We deduce from this description of $A_u(I_2)$ that we have:

\begin{theorem}
The coefficients of the representation $r_{\beta\alpha}=\bar{u}\otimes u-1$ of $A_u(I_2)$ commute, and generate a $C^*$-algebra equal to $C(SO(3))$. A similar result holds for $r_{\alpha\beta}=u\otimes\bar{u}-1$.
\end{theorem}

\begin{proof}
By using Lemma 7 we see that the representation $\bar{u}\otimes u$ of $A_u(I_2)_{red}$ can be indentified with the following representation of $C(\mathbb T)*_{red}C(SU(2))$:
$$\overline{\begin{pmatrix}za&zb\\-z\bar{b}&z\bar{a}\end{pmatrix}}\otimes
\begin{pmatrix}za&zb\\-z\bar{b}&z\bar{a}\end{pmatrix}
=\overline{\begin{pmatrix}a&b\\-\bar{b}&\bar{a}\end{pmatrix}}\otimes
\begin{pmatrix}a&b\\-\bar{b}&\bar{a}\end{pmatrix}$$

Thus the representation $\bar{u}\otimes u-1$ of $A_u(I_2)_{red}$ corresponds to the 3-dimensional representation of $SU(2)$, i.e. to the fundamental representation of $SO(3)$. It follows that the coefficients of $r_{\beta\alpha}=\bar{u}\otimes u-1$ of $A_u(I_2)_{red}$ pairwise commute, and generate a commutative $C^*$-algebra, equal to $C(SO(3))$. We can conclude by using the fact that $C(SO(3))$ is the enveloping $C^*$-algebra of $C(SO(3))_s$.
\end{proof}

We show now that the von Neumann algebra $A_u(I_2)_{red}''$ (the bicommutant of the image of $A_u(I_2)$ by the left regular representation on $l^2(A_u(I_2),h)$, $h$ being the Haar measure) is isomorphic to the factor $W^*(F_2)$ associated to the free group on two generators $F_2$.

We recall that if $(M,\phi )$ is a unital $*$-algebra with a unital linear form and $A\subset M$ is a unital $*$-subalgebra, an element $x\in M$ is called $*$-free from $A$ when the $*$-algebra generated by $x$ inside $M$ is free from $A$. We will need the following technical lemma:

\begin{lemma}
Let $(M,\phi)$ be a $*$-algebra with a trace, $1\in A\subset M$ be a $C^*$-subalgebra, $d\in A$ be a unitary such that $\phi(d)=\phi(d^*)=0$, and $u\in M$ be a Haar unitary which is $*$-free from $A$. Then $ud$ is a Haar unitary, $*$-free from $A$.
\end{lemma}

\begin{proof}
The element $ud$ in the statement is clearly a Haar unitary. In order to prove that $ud$ is $*$-free from $A$ we have to verify that if $a_i\in\mathbb Z^*$ and $f_i\in A\cap\ker(\phi)$
($1\leq i\leq n$), then $P=(ud)^{a_1}f_1(ud)^{a_2}f_2\ldots$ belongs to $ker(\phi)$. 

$P$ is a product of $u,u^*$ terms alternating with $d,d^*,f_i,df_i,f_id^*,df_id^*$ terms. Observe that $\phi (d)=\phi (d^*)=\phi (f_i)=\phi(df_id^*)=0$, and that the terms of the form $df_i$ and $f_id^*$ appear inside $P$ between $u$ and $u$, or between $u^*$ and $u^*$. 

We write each $df_i$ as $[\phi (df_i)1]+[df_i-\phi (df_i)1]$, and each $f_id^*$ as $[\phi (f_id^*)1]+[f_id^*-\phi (f_id^*)]$. By developing $P$, we obtain a certain linear combination of products of elements of the form $u^k$ with $k\in\mathbb Z^*$, alternating with elements of $A\cap \ker(\phi )$. Thus $\phi (P)=0$. 
\end{proof}

We can now prove our von Neumann algebra result:

\begin{theorem}
$A_u(I_2)_{red}''=W^*(F_2)$.
\end{theorem}

\begin{proof}
By Lemma 7 we have an embedding as follows:
$$A_u(I_2)_{red}''\subset L^{\infty}(\mathbb  T)*L^{\infty} (SU(2))\, ,\,\,\,
u\to\begin{pmatrix}za&zb\\-z\bar{b}&z\bar{a}\end{pmatrix}$$ 

Let $d=sgn\circ (a+\bar{a})$, i.e. the composition of $a+\bar{a}:SU(2)\to\mathbb R$ with the sign function $sgn:\mathbb R\to\{-1,0,1\}$. Then $d\in L^{\infty} (SU(2))$ is a unitary such that $d^2=1$.

The polar part of $za+z\bar{a}$ being $zd$, we have $zd\in A_u(I_2)_{red}''$. Thus $dz^*\in A_u(I_2)_{red}''$, and by multiplying at left by $dz^*$ the generators $za,zb,z\bar{a},z\bar{b}$ of $A_u(I_2)_{red}''$, we obtain that  $A_u(I_2)_{red}''$ is generated by $zd$, $da$, $db$, $d\bar{a}$, $d\bar{b}$.

By using Lemma 8, $W^*(da,db,d\bar{a},d\bar{b})$ and $W^*(zd)$ are diffuse abelian algebras which generate $A_u(I_2)_{red}''$, and this implies $A_u(I_2)_{red}''=W^*(F_2)$ (see Theorem 2.6.2 in \cite{vdn}).
\end{proof}

{\bf Remark.} Let $U_n^{nc}$ be the universal $C^*$-algebra generated by the entries of a $n\times n$ unitary matrix. By combining the formula $U_{n,red}^{nc}\otimes M_n=M_n*_{red}C(\mathbb T)$ of McClanahan \cite{mc1} with the formula $M_n*W^*(F_s)=W^*(F_{n^2s})\otimes M_n$ of Dykema \cite{dyk} we obtain $U_{2,red}^{nc,"}=W^*(F_4)$.

\section{Remarks on the adjoint representation}

The adjoint representation of a discrete group $\Gamma$ is an important tool for studying various questions regarding the algebras $C^*(\Gamma)$ and $C^*_{red}(\Gamma)$, namely the amenability of $\Gamma$, the nuclearity of $C^*_{red}(\Gamma)$, the simplicity of $C^*_{red}(\Gamma)$, and so on. We will study here this type of questions for the $C^*$-algebras $A_u(F)$ and $A_u(F)_{red}$, which, from the Pontrjagin duality viewpoint, are the ``quantum $C^*(F_n)$ and $C^*_{red}(F_n)$''.

If $G$ is a compact matrix quantum group {\em whose Haar measure is a trace} we can define on $G_{full}$ the left and right regular representations $\lambda,\rho$ by $\lambda (x)(y)=xy$ and $\rho (x)(y)=y\kappa (x)$, and the adjoint representation as being the following composition:
$$G_p\displaystyle{\mathop{\longrightarrow}^\delta}G_{full}\otimes_{max}G_{full}\displaystyle{\mathop{\longrightarrow}^{\lambda\otimes\rho}}B(l^2(G_{red}))$$

The general case is more subtle, and we will use the following theorem of Woronowicz (see e.g. the Rappels 5.1 in \cite{bsk}):

\begin{theorem}[\cite{wo2}]
To any compact matrix quantum group $G$ we can can associate a (unique) family of characters $(f_z)_{z\in\mathbb C}$ of $G_s$ satisfying:
\begin{enumerate}
\item $h(ab)=h(b (f_1*a*f_1))$, for any $a\in G_s$ and $b\in G$.

\item $\kappa^2(a)=f_{-1}*a*f_1$, for any $a\in G_s$.

\item $f_0=e$ (the counit of $G_s$) and $f_{z+t}=f_z*f_t$, for any $z,t\in\mathbb C$.

\item $f_z*\kappa (a)=\kappa (a*f_{-z})$ and $\kappa (a)*f_z=\kappa (f_{-z}*a)$ for any $a\in G_s$ and $z\in\mathbb C$.
\end{enumerate}
In addition, these characters $f_z$ can be defined as follows: if $u\in\widehat{G}$ and $F$ is the unique positive matrix intertwining $u$ and $(I\otimes\kappa )(u^*)=(I\otimes\kappa^2)(u)$, normalised as to satisfy $Tr(F)=Tr(F^{-1})$ (cf. Theorem 5.4 in \cite{wo2}), then $(Id\otimes f_z)(u)=F^z$.
\end{theorem}

{\bf Notations.} $\mathcal L(G_{red})$ is the algebra of bounded operators $G_{red}\to G_{red}$. The character map $\widehat{G}\to
G_s$ will be denoted $\chi$ (see Reminders (B) and (C)).

\medskip

As a consequence of the above result, we have:

\begin{corollary}
Let $G$ be a compact matrix quantum group, and $a,b,c,d\in\mathbb C$.
\begin{enumerate}
\item $x\mapsto f_a*x*f_b$ is an automorphism of $G_s$, and $\lambda_{a,b}:G_s\to\mathcal L(G_{red})$ given by $\lambda_{a,b}(x)(y)=(f_a*x*f_b)y$ is a morphism of unital $C^*$-algebras.

\item $x\to f_c*\kappa (x)*f_d$ is an antiautomorphism of $G_s$, and 
$\rho_{c,d}:G_s\to\mathcal L(G_{red})$ given by $\rho_{c,d}(x)(y)=y(f_c*\kappa (x)*f_d)$ is a morphism of unital $C^*$-algebras.

\item The map $G_s\displaystyle{\mathop{\longrightarrow}^\delta}G_s\otimes G_s
\displaystyle{\mathop{\longrightarrow}^{\lambda_{a,b}\otimes\rho_{c,d}}} 
{\mathcal L}(G_{red})$ is a morphism of unital $C^*$-algebras.
\end{enumerate}
\end{corollary}

The morphism defined in (3) allows to associate to any 
$a,b,c,d\in\mathbb C$ a certain map $ad:\widehat{G}\to\mathcal L(G_{red})$. The interest of the following lemma will appear later on, in the simplicity proofs, where we will use several maps $ad$ of this type:

\begin{lemma}
Let $G$ be a compact matrix quantum group, let $a,b,c,d\in\mathbb C$, and set:
$$ad=(\lambda_{a,b}\otimes\rho_{c,d})\circ \delta\circ\chi :\widehat{G}\to
{\mathcal L}(G_{red})$$
Let also $u\in\widehat{G}$ and let $F$ be the associated matrix, from Theorem 7. Then: 
\begin{enumerate}
\item $ad(u)(z)=\sum_{i,k}(F^buF^a)_{ik}z(F^{-c}u^*F^{-d})_{ki}$.

\item If $a+c=b+d=0$, then $ad(u)$ is of type $z\to\sum a_kza_k^*$ with $a_k\in G_s$ (finite sum).

\item If $a=c$, then $ad(u)(1)=K\cdot 1$ with $K\in\mathbb R_+^*$.

\item If $s,t\in\mathbb R$ satisfy $t+b-d=0$ or $s+a-c+1=0$, there exists $M\in\mathbb R_+^*$ such that $ad(u)/M$ preserves any state $\phi$ of $G_{red}$ satisfying $\phi (xy)=\phi (y(f_s *x*f_t))$. 
\end{enumerate}
\end{lemma}

\begin{proof}
(1) By using Theorem 7, we have:
$$f_a*u_{ij}*f_b=(f_b\otimes Id\otimes f_a)(\sum u_{is}\otimes u_{sk}\otimes u_{kj})=(F^buF^a)_{ij}$$

By Theorem 7 (4), and by using  $(I\otimes\kappa )(u)=u^*$, we obtain:
$$f_c*\kappa (u_{ij})*f_d=\kappa (f_{-d}*u_{ij}*f_{-c})=(F^{-c}u^*F^{-d})_{ij}$$

Thus $(\lambda_{a,b}\otimes\rho_{c,d})(u_{ik}\otimes u_{lj})(z)=(F^buF^a)_{ik}z(F^{-c}u^*F^{-d})_{lj}$, and so:
$$(\lambda_{a,b}\otimes\rho_{c,d})\delta (u_{ij})(z)=\sum_k (F^buF^a)_{ik}z(F^{-c}u^*F^{-d})_{kj}$$

(2) We have $(F^{-c}u^*F^{-d})_{ki}=[(F^{-c}u^*F^{-d})^*]_{ik}^*=(F^{-d}uF^{-c})_{ik}^*$ (recall that $F>0$), and so if $a+c=b+d=0$, then:
$$ad(u):z\to \sum_{i,k} (F^buF^a)_{ik}z(F^buF^a)_{ik}^*$$

(3) If $a=c$ then $ad(u)(1)=\sum_{i,k}(F^buF^a)_{ik}1(F^{-a}u^*F^{-d})_{ki}=Tr(F^{b-d})1$.

(4) We have $\phi (ad(u)(z))=\phi [\sum (f_a*u_{ik}*f_b)z(f_c*\kappa (u_{ki})*f_d)]=\phi (z M)$, where:
$$M=\sum_{i,k}(f_c*\kappa (u_{ki})*f_d)(f_{s +a}*u_{ik}*f_{t+b})$$

If $t+b-d=0$, by using the formulae in the proof of (1) we obtain:
$$M=\sum_{i,k}(F^{-c}u^*F^{-d})_{ki}(F^{t +b}uF^{s +a})_{ik}=Tr(F^{s +a-c})>0$$

Assume now $s +a-c+1=0$. By using Theorem 7 (2) we have $u_{ik}=f_1*\kappa^2(u_{ik})*f_{-1}$. Together with Theorem 7 (4) and with the formulae in the proof of (1), this gives:
\begin{eqnarray*}
M
&=&\sum_{i,k}(f_c*\kappa (u_{ki})*f_d)(f_{s+a+1}*\kappa^2(u_{ik})*f_{t +b-1})\\
&=&\kappa [\sum_{i,k}(f_{-t -b+1}*\kappa (u_{ik})*f_{-s
-a-1})(f_{-d}*u_{ki}*f_{-c})]\\
&=&\kappa [(I\otimes Tr)(F^{t +b-1}u^*F^{s +a+1-c}uF^{-d})]\\
&=&Tr(F^{t+b-1-d})\\
&>&0
\end{eqnarray*}

Thus $M$ is positive in both cases, and we are done.
\end{proof}

\section{Powers' property}

Let $(A,\tau)$ be a unital $C^*$-algebra with a faithful trace. Haagerup et Zsido proved in \cite{hzs} that $A$ is simple with unique trace if and only if it has the Dixmier property:

\medskip

{\em For any $a\in A$, the closed convex envelope of $\{uau^*|u\ {\rm unitary}\}$ contains a scalar multiple of the unit $1_A$}.

\medskip

The simplicity of $C^*_{red}(F_n)$ was proved in \cite{pow}, and Powers' method has been extended to various free products \cite{avi}, \cite{mc2}, or $C^*$-algebras of discrete groups \cite{bch}, \cite{har}, \cite{hsk}. The proofs use technical estimates inside $B(l^2(A,\tau ))$, which ``move towards 0'' any element having trace 0, by using sums of inner automorphisms, i.e. which prove the Dixmier property.

In the case of the reduced $C^*$-algebras of discrete groups $\Gamma$, the estimate inside $l^2(\Gamma )$ is usually obtained by using combinatorial or geometric properties of $\Gamma$. One of such properties is Powers' property, defined in \cite{har}:

\medskip

{\em For any finite subset $F\subset\Gamma -\{ 1\}$, there exist elements
$g_1,g_2,g_3\in\Gamma$ and a partition $\Gamma=D\coprod E$ such that $F\cdot D\cap D=\emptyset$ and $g_s\cdot E\cap g_k\cdot E=\emptyset,\forall\, s\neq k$}.

\medskip

It is easy to see that the free groups $F_n$ have Powers' property. This property appears in many other contexts (see \cite{har}), for instance any action of $\Gamma$ on a Hausdorff space, which is strongly hyperbolic, minimal and strongly faithful, produces a partition $\Gamma =D\coprod E$ and elements (in fact, an infinity of elements) $g_i$ as above. Under this assumption, the proof of simplicity of $C^*_{red}(\Gamma)$ from \cite{hsk} has two steps, as follows:
\begin{enumerate}
\item If $x\in l^2(F)$ is self-adjoint of trace 0, then $||1/3\sum
u_{g_s}xu_{g_s}^*||\leq 0.98||x||$.

\item $C^*_{red}(\Gamma)$ has the Dixmier property, and is therefore simple, with unique trace.
\end{enumerate}

We will extend this proof to the following compact quantum groups:

\begin{definition}
Let $G$ be a compact quantum group. We endow the set $P(\widehat{G})$ of parts of $\widehat{G}$ with the involution $\bar{A}=\{\bar{a}|a\in A\}$ and with the multiplication $\circ$ defined by:
$$A\circ B=\{ r\in\widehat{G}\mid\exists\, a\in A,\exists\, b\in B\,\,\, {\rm with}\,\,\, r\subset a\otimes b\}$$
We say that $G$ has Powers' property if for any $F\subset\widehat{G}-\{ 1\}$ finite there exist $r_1,r_2,r_3\in\widehat{G}$ and a partition $\widehat{G}=D\coprod E$ such that $F\circ D\cap D=\emptyset$ and $r_s\circ E\cap r_k\circ E=\emptyset$, $\forall\, s\neq k$.
\end{definition}

To be more precise, the aim of this last part of the paper is to prove that the algebras $A_u(F)_{red}$ are simple (with at most one trace). We have to solve three questions:
\begin{enumerate}
\item Extend the simplicity proof from \cite{hsk} to the compact quantum groups having Powers' property, and a tracial Haar state.

\item Extend (1) to the arbitrary compact quantum groups having Powers' property.

\item Extend (2) to $A_u(F)$, which does not have Powers' property. Indeed, Theorem 1 shows that $F=\{ r_{\alpha\beta},r_{\beta\alpha}\}$ satisfies $r_x\in F\circ \{ r_x\}$, for any $1\neq r_x\in \widehat{A_u(F)}$.
\end{enumerate}

In fact, all our proofs of simplicity and of non-existence of KMS states will be based on the same estimate (Proposition 8 below). Observe that:

-- If the Haar measure of $G$ is a trace, the statement of Proposition 8 considerably simplifies. It goes the same for the proof: we do not need the material from section 6.

-- If in addition we have $G=C^*_{red}(\Gamma )$ (e.g. $C^*_{red}(F_2)$), then Proposition 8 is the technical lemma in \cite{hsk}, but the estimate that we obtain here is stronger, namely:
$$||\frac{1}{3}\sum_{s=1,2,3} u_{g_s}xu_{g_s}^*||\leq \frac{2\sqrt{2}}{3}||x||$$

Finally, observe that the methods that we develop here have no chances to apply to $A_o(F)$, because the algebra $A_o(F)_{central}$ is commutative. 

\medskip

{\bf Notation.} For $x\in G_s$ we denote by $supp(x)\subset P(\widehat{G})$ the set of irreducible representations having coefficients appearing in $x$ (we recall that the space of coefficients of $r\in B(H_r)\otimes G_s$ is $G(r)= \{ (\phi\otimes Id)(r)|\phi\in B(H_r)^*\}$; by \cite{wo2} we have $G_s=\oplus_{r\in\widehat{G}} G(r)$).

\medskip

We prove now our main estimate:

\begin{proposition}
Let $(G,u)$ be a reduced compact matrix quantum group, and let $s,t\in\mathbb R$. Let $\widehat{G}=D\coprod E$ be a partition, and $r_1,r_2,r_3\in\widehat{G}$ be such that $r_l\circ E\cap r_k\circ E=\emptyset$, $\forall\, l\neq k$. Then there exists a unital linear map  $T:G\to G$ such that:
\begin{enumerate}
\item There exists a finite family $\{ a_i\}$ of elements of $G_s$ such that $T:z\to \sum a_iza_i^*$.

\item $T$ preserves the states $\phi\in G_{red}^*$ satisfying $\phi
(xy)=\phi (y(f_s *x*f_t )),\,\forall\,  x,y\in G_s$.

\item For any $z=z^*\in G_s$ such that $supp(z)\circ D\cap D=\emptyset$, we have $||T(z)||\leq0.95||z||$, and $supp(T(z))\subset \cup_i r_i\circ supp(z)\circ \bar{r}_i$.
\end{enumerate}
\end{proposition}

\begin{proof}
Lemma 9 with $a=c=0$ and $d=-b=t /2$ gives a certain map $ad:\widehat{G}\to\mathcal L(G)$. Our choice of $a,b,c,d$ implies $a+c=b+d=0$, $a=c$ and $t +b-d=0$, so we can use Lemma 9 (2,3,4) with $u=r_i$, for $i=1,2,3$. We obtain three finite families $\{ a_{i,k}\}_k$ of elements of $G_s$ and six positive numbers $K_i,M_i$ ($i=1,2,3$) such that, for any $i$:

(i) $ad(r_i)(z)=\sum_k a_{i,k}za_{i,k}^*$.

(ii) $ad(r_i)(1)=K_i\cdot 1$.

(iii) $ad(r_i)/M_i$ preserves $\phi$.

Since $\phi (ad(r_i)(1))=M_i$ (by (iii)) and $\phi (ad(r_i)(1))=K_i$ (by (ii)), we have $K_i=M_i$ for any $i$. We construct now a linear map $T$, as follows: 
$$T=\frac{1}{3}\left( \frac{ad(r_1)}{M_1}+\frac{ad(r_2)}{M_2}+\frac{ad(r_3)}{M_3}\right)$$

Since (1,2) are satisfied, it remains to prove $(3)$. Let $h$ be the Haar measure of $G$, and $(H,\pi
)$ be the GNS construction for $(G,h)$. For $i=1,2,3$ we set:
$$T_i':B(H)\to B(H),\,\, P\to M_i^{-1}\sum_k \pi (a_{i,k})P\pi
(a_{i,k}^*)$$

Let $T'=(T_1'+T_2'+T_3')/3$. This is a completely positive unital map. Let $p,q$ be the projections on $H$ onto the closure of the linear spaces generated by the coefficients of the representations in $D,E$. We have $\widehat{G}=D\coprod E$ and $supp(z)\circ D\cap D=\emptyset$, and so:
$$p+q=1\, ,\, p\pi (z)p=0$$

If $t\neq s\in \{1,2,3\}$, then $\bar{r}_t\circ r_s\circ E\cap E=\emptyset$: indeed, if $r,p\in E$ are such that $r\subset\bar{r}_t\otimes r_s\otimes p$, then $h(\chi (\overline{r_t\otimes r}\otimes r_s\otimes p))\geq 1$, so
$r_t\otimes r$ and $r_s\otimes p$ have a common irreducible component, which must be in $r_s\circ E\cap r_t\circ E=\emptyset$, contradiction (cf. Reminders (C)). 

Thus if $a,b$ are arbitrary coefficients of  $r_t,r_s$, then $q\pi (a^*b)q=0$. But the elements $a_{i,k}$ being coefficients of $r_i$ for any $i$ (cf. Lemma 9 (1)), we have:
$$T_t'(q)\, T_s'(q)=(M_tM_s)^{-1}\sum_{k,h}\pi (a_{t,k})q\pi
(a_{t,k}^*a_{s,h})q\pi (a_{s,h}^*)=0$$ 

By using this, and the fact that $(T_i')^n$ are completely positive and unital, we obtain:
$$||T'(q)||=\lim_{n\to\infty}||T^{\prime}(q)^n||^{\frac{1}{n}}=
\frac{1}{3}\lim_{n\to\infty}||(\sum_i T_i^{\prime}(q))^n||^{\frac{1}{n}}=
\frac{1}{3}\lim_{n\to\infty}||\sum_i (T_i^{\prime}(q))^n||^{\frac{1}{n}}
\leq\frac{1}{3}$$

The assertion (3) on $supp(T(z))$ is clear. Finally, the estimate $||T(z)||\leq 0.95||z||$ follows from the next lemma (applied with $f=T',\, x=\pi (z)$ and $\delta =1/3$) and from the fact that the GNS representation $\pi$ is isometric.
\end{proof}

\begin{lemma}
Let $H$ be a Hilbert space, $x=x^*\in B(H)$, $p+q=1$ projections on $H$, and $f:B(H)\to B(H)$ be completely positive and unital. If $pxp=0$
and $||f(q)||\leq\delta <1/2$, then $||f(x)||\leq 2\sqrt{\delta -\delta^2}||x||$.
\end{lemma}

\begin{proof} (G. Skandalis) Let $\zeta\in H$ be arbitrary, of norm 1. We want to prove that we have $|<f(x)\zeta ,\zeta >|\leq\, 2\sqrt{\delta -\delta^2}|x|$. By using the  Stinespring theorem we can assume $f(z)=\omega ^*z\omega$ with $\omega ^*\omega=1$. With $\xi =\omega\zeta$, it is enough to show that:

\medskip

{\em If $x=x^*\in B(H)$, $p+q=1$ are projections on $H$, with $pxp=0$, and $\xi\in H$ has norm $1$ and satisfies $<q\xi ,\xi >\,\leq\, \delta <1/2$, then $|<x\xi ,\xi >|\leq\, 2\sqrt{\delta -\delta^2}||x||$.}

\medskip

Let $E\in B(H)$ be the projection onto $\mathbb C p\xi\oplus \mathbb C q\xi$. We can replace in the above statement $H,p,q,x,\xi$ by$E(H),EpE,EqE,ExE,\xi$. Indeed, we have $<q\xi ,\xi>=<EqE\xi ,\xi >$ and $<x\xi ,\xi >=<ExE\xi,\xi >$ and $||ExE||\leq||x||$. We can assume as well $||x||=1$. 

Thus, assume $H=\mathbb C^2$, $p=\begin{pmatrix}0&0\\0&1\end{pmatrix}$, $q=\begin{pmatrix}1&0\\0&0\end{pmatrix}$, $x=\begin{pmatrix}a&b\\ \bar{b}&0\end{pmatrix}$ and $\xi=\begin{pmatrix}m\\ n\end{pmatrix}$ with $a\in\mathbb R$ and $b,m,n\in\mathbb C$. We have:
$$<x\xi ,\xi >=\left<\begin{pmatrix}am+bn\\ \bar{b}m\end{pmatrix},\begin{pmatrix}m\\ n\end{pmatrix}\right>=
a|m|^2+2Re(b\, n\,\bar{m})$$

We have $|m|^2=<q\xi ,\xi >\,\leq\, \delta$ and $|m|^2+|n|^2=||\xi||^2 =1$, and so:
$$|<x\xi ,\xi >|\leq\delta|a|+2\sqrt{\delta(1-\delta )}|b|$$

We can assume $a\geq 0$. The roots of $\det(x-zI)=z^2-az-|b|^2$ are
$(a\pm\sqrt{a^2+4|b|^2})/2$. But $||x||=1$, so these roots are in $[-1,1]$, so $\sqrt{a^2+4|b|^2}\leq 2-a$, so $a\leq1-|b|^2$. Thus:
$$|<x\xi ,\xi >|\leq\delta (1-|b|^2) +2\sqrt{\delta (1-\delta
)}|b|=1-(\sqrt{1-\delta}-\sqrt{\delta}|b|)^2$$

We have $\delta <1/2$, so the function $b\to1-(\sqrt{1-\delta}-\sqrt{\delta}|b|)^2$ reaches its maximum on $[-1,1]$ at $b=\pm 1$. But this maximum is $2\sqrt{\delta -\delta^2}$, and we are done.
\end{proof}

Finally, we will use the following technical lemma, instead of the Dixmier property:

\begin{lemma}
Let $(A,\phi )$ be a $C^*$-algebra with a faithful state, let $\psi\in A^*$ be a state, let $A_s\subset A$ be a dense $*$-subalgebra, and $0<\delta <1$. If for any $x\in ker(\phi )\cap A_s$ self-adjoint there exists a finite family of elements $a_i\in A$ such that $z\to\sum a_iza_i^*$ is unital, preserves $\phi$ and $\psi$, and maps $x$ to an element of norm $\leq\delta ||x||$, then $A$ is simple, and $\psi =\phi$.
\end{lemma}

\begin{proof}
We can assume $A_s=A$. By applying several times the condition in the statement, we can assume that $\delta >0$ is as small as we want.

Let $J\subset A$ be a two-sided ideal, $y\in J$ and $z=yy^*/\phi (yy^*)$. Then we can find elements $a_i$ with $||\sum a_i(1-z)a_i^*||< 1$, i.e. with $\sum a_iza_i^*$ invertible. But $\sum a_iza_i^*\in J$, so $J=A$.

Let $x=x^*\in\ker(\phi )$ arbitrary, and $\epsilon >0$ small. We can find an element $y=\sum a_ixa_i^*$ having norm $<\epsilon$, so $|\psi (x)|=|\psi (y)|\leq\epsilon$. Thus $\psi (x)=0$, and so $\psi =\phi$ on the self-adjoints. Now since any operator is a finite linear combination of 1 and of self-adjoints belonging to $\ker(\phi )$, we have $\psi =\phi$.
\end{proof}

We can go back now to Powers' property, as defined above, and prove:

\begin{proposition}
If $G$ has Powers' property then $G_{red}$ est simple.

Assume in addition that we are given a state $\psi$ of $G_{red}$ such that $\forall\, x,y\in G_s$ we have $\psi (xy)=\psi (y(f_1 *x*f_1 ))$. Then $\psi =h$ (the Haar measure of $G_{red}$).

This if $h$ is a trace, t is the unique trace of $G_{red}$. 
\end{proposition}

\begin{proof}
Let $x\in\ker(h)\cap G_s$ be self-adjoint. Observe that $1$ is not in  $F=supp(x)$. By using Powers' property, we can apply Proposition 8 with $s =t =1$. We obtain in this way a unital map $f$ of the form $z\to \sum a_i za_i^*$ which leaves invariant both $h$ and $\psi$ (cf. Proposition 8 (2)), and which is such that $||f(x)||\leq 0.95||x||$ (norms inside $G_{red}$). We can conclude by using Lemma 11 (with $A=G_{red},\, A_s=G_s$ and $\phi =h$).
\end{proof}

\section{Simplicity of $A_u(F)_{red}$}

$A_u(F)_{red}$ does not have Powers' property (take indeed $F=\{
r_{\alpha\beta},r_{\beta\alpha}\}$), but we will show however that this algebra is simple, by using Proposition 8. For this purpose, we need to identify the objects used in section 7, for $G=A_u(F)$. By using our description of the representations of $A_u(F)$, we can identify $\widehat{A_u(F)}=\mathbb N*\mathbb N$. The
multiplication $\circ$ on $P(\mathbb N*\mathbb N)$ is given (for
$x,y\in\mathbb N*\mathbb N $) by the following formula:
$$\{ x\}\circ \{ y\} =\{ ab\mid\exists\, g\in \mathbb N *\mathbb N \,\,\, {\rm with}\,\,\, x=ag,\,y=\overline{g}b\}$$

{\bf Notation.} For $w\in\mathbb N*\mathbb N$ we denote by $\{w\ldots\}$ (resp. $\{\ldots w\}$) the set of words in $\mathbb N*\mathbb N$ which begin (resp. end) with $w$. For $w,y\in\mathbb N*\mathbb N$ we set $\{ w\ldots y\}=\{ w\ldots\}\cap \{\ldots y\}$. We denote by $(\beta\alpha)^N$ the word $\beta\alpha\beta\alpha\ldots\beta\alpha$ ($N$ times).

\medskip

The following lemma is clear:

\begin{lemma}
Consider the sets $D=\{\alpha\ldots\}$, $E=\{\beta\ldots\}\cup\{e\}$, $F=\{ \beta\ldots\alpha\}$ and the elements $r_1=\beta\alpha\beta$, $r_2=\beta\alpha^2\beta$ and $r_3=\beta\alpha^3\beta$. Then $\mathbb N*\mathbb N  =D\coprod E$ is a partition, $F\circ D\cap D=\emptyset$, and $r_s\circ E\cap r_k\circ E=\emptyset$, $\forall\, s\neq k$.
\end{lemma}

We fix $n\in\mathbb N$ and $F\in GL(n,\mathbb C)$, and we set $G=A_u(F)_{red}$, with Haar measure denoted $h$. The estimate in Proposition 8 above has the following consequence:

\begin{corollary}
Let $s,t\in\mathbb R$, $\epsilon >0$, and $\psi$ be a state of $G$ such that $\psi (xy)=\psi(y(f_s *x*f_t ))$, $\forall\, x,y\in G_s$. Then there exists a unital linear map $V:G\to G$ of the form $z\to\sum a_iza_i^*$ (finite sum, with $a_i\in G_s$) which preserves $\psi$, such that for any $x=x^*\in A_u(F)_s$ with $supp(x)\subset \{ \beta\ldots\alpha\}$ we have $||V(x)||\leq\epsilon||x||$ and $supp(V(x))\subset \{ \beta\ldots\alpha\}$.
\end{corollary}

\begin{proof}
We apply Proposition 8 to the parts constructed in Lemma 12. We obtain a certain map $T:G\to G$ of the form $z\to\sum b_izb_i^*$ (finite sum, with $b_i\in G_s$) which preserves $\psi$. By setting $z=x$ in Proposition 8 (3) we obtain $||T(x)||\leq 0.95||x||$ and:
$$supp(T(x))\subset \bigcup_i r_i\circ supp(x)\circ \overline{r_i}\subset\{ \beta
\ldots\beta\}\circ \{\beta\ldots\alpha \}\circ \{ \alpha\ldots\alpha \}\subset \{\beta\ldots\alpha \}$$

By Proposition 8 (1) we have $T(x)=T(x)^*\in G_s$, so we can apply Proposition 8 (3) with $z=T(x)$, then with $z=T^2(x)$, and so on. Thus, in order to finish, we can choose $m\in\mathbb N$ such that $0.95^m\leq\epsilon$, and set $V=T^m$. 
\end{proof}

We will need the following combinatorial lemma:

\begin{lemma}
$\forall\, F\subset \mathbb N *\mathbb N $ finite, $(\beta\alpha )^N\circ F\circ
(\beta\alpha )^N\subset\{\beta\ldots\alpha\}\cup\{ e\}$ for $N$ big enough.
\end{lemma}

\begin{proof}
If we denote by $Y\subset\mathbb N*\mathbb N$ the set of alternating words (i.e. words which do not contain $\alpha^2$ or $\beta^2$), it is easy to see that:

$(a)$ $ Y \circ\{ ...\alpha\}\cap\{ ...\beta\}=\emptyset$ ; $(b)$ $ Y \circ\{ ...\beta\}\cap\{
...\alpha\}=\emptyset$

$(c)$ $\{\alpha ...\}\circ Y \cap\{\beta ...\}=\emptyset$ ; $(d)$ $\{\beta
...\}\circ Y \cap\{\alpha ...\}=\emptyset$. 

It suffices to prove the lemma for 1-element sets, $F=\{ z\}$, and here we have:

-- Assume $z\in Y$. By $(d)$, we have that $(\beta\alpha )^N\circ z$ equals $e$ (and in this case, we are done), or starts with $\beta$. By using one more time $(d)$ we see that $(\beta\alpha )^N\circ z\circ (\beta\alpha )^N$ equals $e$ or begins with $\beta$. Similarly, by using twice $(a)$, we see that $(\beta\alpha
)^N\circ z\circ (\beta\alpha )^N$ equals $e$ or ends with $\alpha$. Thus $(\beta\alpha )^N\circ z\circ (\beta\alpha )^N$ belongs to $\{\beta\ldots\alpha\}\cup\{e\}$.

-- Assume $z\in \mathbb N*\mathbb N-Y$, for instance $z=x\alpha^2y$. Then $(\beta\alpha )^N\circ x\alpha\subset \{\ldots\alpha\}\cup\{ e\}$ by $(a)$. For $N\geq l(x)$, it is clear that $(\beta\alpha )^N\circ x\alpha\subset\{\beta 
\ldots\alpha\}$. By the same arguments, $\alpha y\circ (\beta\alpha )^N\subset\{\alpha\ldots\alpha\}$ for $N\geq l(y)$. Thus for $N$ big enough:
$$(\beta\alpha )^N\circ (x\alpha^2y)\circ (\beta\alpha )^N=[(\beta\alpha )^N\circ x\alpha ]\circ [\alpha y\circ (\beta\alpha )^N]\subset\{\beta
...\alpha\} $$

Summarizing, the conclusion holds in both cases, and we are done.
\end{proof}

The above results have the following consequence:

\begin{corollary}
Let $x=x^*\in G_s$ be such that $h(x)=0$, and let $v,w\in\mathbb R$.
\begin{enumerate}
\item There exists a unital linear map $W:G\to G$ of the form $z\to\sum b_izb_i^*$ (finite sum, with $b_i\in G_s$), which preserves $h$, and which satisfies $supp(W(x))\subset \{\beta\ldots\alpha\}$.

\item There exist $L\in\mathbb R_+^*$ and $U:G\to G$ linear, of the form $z\to\sum c_izc_i^*$ (finite sum, with $c_i\in G_s$), which preserves $h$, such that $supp(U(x))\subset\{\beta\ldots\alpha\}$ and such that $U/L$ preserves any state $\psi$ of $G$ satisfying $\psi (pq)=\psi (q(f_v*p*f_w))$, $\forall\, p,q\in G_s$. 
\end{enumerate}
\end{corollary}

\begin{proof}
Pick $K\in\mathbb N$, $(\beta\alpha )^K\circ supp(x)\circ
(\beta\alpha )^K\subset \{\beta\ldots\alpha\}\cup\{ e\}$, and set $r=r_{(\beta\alpha )^K}$.

(1) Lemma 9 with $a=c=0$ and $d=-b=1/2$ gives a certain map $ad :\widehat{G}\to\mathcal L(G)$. This choice of $a,b,c,d$ allows us to apply (with $u=r$) Lemma 9 (2,3), as well as Lemma 9 (4), with $s=t=1$ and $\phi =h$. We obtain two numbers $K,M>0$, which are equal.

With $W=ad(r)/M$, it remains to show that $supp(W(x))\subset \{\beta\ldots\alpha\}$. By using $r\circ supp(x)\circ r\subset \{\beta\ldots\alpha\}\cup\{ e\}$, the formula of $ad$ (Lemma 9 (1)), and $\overline{r}=r$ we obtain $supp(W(x))\subset\{\beta\ldots\alpha\}\cup\{ e\}$. But $h(W(x))=h(x)=0$, so $e\not\in supp(W(x))$.

(2) Lemma 9 applied with $c=-a=(v+1)/2$ and $d=-b=1/2$ gives a certain map  $ad:\widehat{G}\to\mathcal L(G)$. We can apply (with $u=r$) Lemma 9 (2), as well as Lemma 9 (4), with $s=t=1$ and $\phi=h$. We obtain in this way a certain number $M\in\mathbb R_+^*$ such that $U=ad(r)/M$ is of the form $z\to\sum c_izc_i^*$, and preserves $h$.

Now let us apply again Lemma 9 (4), with the same values of $a,b,c,d$, but this time with $s =v,\, t=w$ and $\phi =\psi$. We obtain in this way a certain number $M_1\in\mathbb R_+^*$ such that $ad(r)/M_1$ preserves $\psi$. We can therefore set $L=M_1/M$. 

Finally, the assertion regarding $supp(U(x))$ can be proved as in (1) above.
\end{proof}

We can now prove our third and last main result, as follows:

\begin{proof}(of Theorem 3) We recall that we have fixed $F\in GL(n,\mathbb C)$, and we have set $G=A_u(F)_{red}$. Let $x=x^*\in G_s$ be such that $h(x)=0$, and let $\epsilon_1 >0$.

(1) By using Corollary 4 (1) with the above $x$ we obtain a certain map $W:G\to G$. By using Corollary 3 with $s=\beta=1$, $\psi=h$ and $\epsilon=\epsilon_1$ we obtain a certain map $V:G\to G$. Observe that $VW$ has the following properties:

-- $VW$ is unital, of the form $z\to \sum_s a_sza_s^*$ (finite sum).

-- $VW$ preserves $h$.

-- $||(VW)(x)||\leq\epsilon_1||x||$.

Lemma 11 (with $A=G$, $A_s=G_s$ and $\psi =\phi =h$) shows then that $G$ is simple.

(2) Let $s,t\in\mathbb R$ and let $\phi$ be a state of $G$ satisfying $\phi
(pq)=\phi (q(f_s*p*f_t))$. By using Corollary 4 (2) with $v=s$, $w=t$ and $\psi =\phi$ we obtain a map $U:G\to G$ and a number $L\in\mathbb R_+^*$. By using Corollary 3 with $\psi=\phi$ and with $\epsilon >0$ such that $\epsilon||U(x)||<L\epsilon_1$ we obtain a certain map $V:G\to G$. Observe that $VU/L$ has the following properties:

-- $VU/L$ preserves $\phi$.

-- $||(VU/L)(x)||\leq\epsilon L^{-1}||U(x)||\leq \epsilon_1$.

By using the first property, and then by letting $\epsilon_1\to 0$ in the second property, we obtain $\phi(x)=0$. Thus $\phi =h$ on the self-adjoints, and so $\phi =h$.

(3) By Theorem 7, the Haar measure of $A_u(F)$ is a trace if and only if $FF^*\in\mathbb C1$.
\end{proof}

Let us investigate now nuclearity questions. We will use:

\begin{proposition}
Let $(G,u)$ be a compact quantum group whose Haar measure is a trace. Then $(G,u)$ is amenable if and only if $G_{red}$ is nuclear.
\end{proposition}

\begin{proof}
Let $J$ be the kernel of the projection $\pi:G_p\to G_{red}$. We analyse the extensions to $G_{full}$ and $G_{red}$ of the above maps $\lambda\otimes\rho$, $\delta$, $e$, defined on $G_s$:

(1) The left-right representation $\lambda\otimes\rho :G_s\otimes
G_s\to B(l^2(G_{red}))$ is a $*$-morphism, which therefore extends into a map $(\lambda\otimes\rho )_{full}:G_{full}\otimes_{max}G_{full}\to B(l^2(G_{red}))$ (see section 6). The kernel of the projection $\pi\otimes
I:G_{full}\otimes_{max}G_{full}\to G_{red}\otimes_{max}G_{full}$ being 
$J\otimes_{max}G_{full}$ (see \cite{was}), we see that $(\lambda\otimes\rho )_{full}$ factorizes via $\pi\otimes I$ into a map $(\lambda\otimes\rho )_{red}:G_{red}\otimes_{max}G_{full}\to B(l^2(G_{red}))$.

(2) The comultiplication $\delta:G_s\to G_s\otimes G_s$ is a $*$-morphism which extends to $G_{full}$ into a map $\delta_{full}:G_{full}\to
G_{full}\otimes_{max}G_{full}$. By composing with the canonical projection
$G_{full}\otimes_{max}G_{full}\to G_{red}\otimes_{max}G_{full}$ we obtain a map $\delta_1:G_{full}\to G_{red}\otimes_{max}G_{full}$.

(3) The comultiplication $\delta_{red}:G_{red}\to G_{red}\otimes_{min}G_{red}$ lifts into a map $\delta_2:G_{red}\to G_{red}\otimes_{min}G_{full}$ (see  Corollaire A.6 in \cite{bsk}).

(4) The counit $e:G_s\to\mathbb C$ extends into a map $e_{full}:G_{full}\to\mathbb C$. 

If we denote by $\tau:T\to<T1,1>$ the canonical state of $B(l^2(G_{red}))$, we obtain in this way a commutative diagram, as follows:
$$\begin{matrix}
G_{red}&\ \displaystyle{\mathop{\longleftarrow}^\pi}&G_{full}&\
\displaystyle{\mathop{\longrightarrow}^{e_{full}}}&\mathbb C\\
&\ \\
\delta_2\downarrow &\ &\ \delta_1\downarrow&\ &\
\uparrow\tau &\\
&\ \\
G_{red}\otimes_{min}G_{full}&\
\displaystyle{\mathop{\longleftarrow}}&G_{red}\otimes_{max}G_{full}&\
\displaystyle{\mathop{\longrightarrow}^{(\lambda\otimes\rho )_{red}
}}&B(l^2(G_{red}))
\end{matrix}$$

Indeed, the commutation on the left follows from the construction of $\delta_1$, $\delta_2$, and the commutation at right can be checked on the generators $u_{ij}$.

It follows that if $G_{red}$ is nuclear, then $\ker(\pi)\subset\ker(e_{full})$, so $G$ is amenable (cf. Proposition 5.5 in \cite{bla}). For the converse, see Remarques A.13 in \cite{bsk}. 
\end{proof}

By applying this result, we obtain:

\begin{corollary}
$A_u(I_n)_{red}$ is not nuclear. The same holds for $A_o(I_n)_{red}$ with $n\geq 3$.
\end{corollary}

{\bf Remark.} Proposition 6 shows that $A_u(F)\sim_{sim}A_u(F')$, with $F'$ diagonal. The universal property of $A_u(F')$ gives a surjection $A_u(F')\to C^*(F_n)$, and by composing with the similarity map, we obtain a surjection $A_u(F)\to C^*(F_n)$. The kernel of this surjection $A_u(F)\to C^*(F_n)$ being a non-trivial ideal of $A_u(F)$, we have $A_u(F)\neq A_u(F)_{red}$.

\end{document}